\documentclass[12pt,leqno]{article}
\usepackage{amsmath,amssymb,amscd}
\numberwithin{equation}{section}

\newtheorem{lemma}{Lemma}[section]
\newtheorem{theorem}{Theorem}[section]
\newtheorem{corollary}{Corollary}[section]

\newcommand\qed{{\unskip\nobreak\hfil\penalty50\hskip2em\vadjust{}
    \nobreak\hfil$\Box$\parfillskip=0pt\finalhyphendemerits=0\par}}

\newcommand{\x}{\times}
\renewcommand{\a}{\alpha}
\renewcommand{\b}{\beta}
\renewcommand{\d}{\delta}
\newcommand{\D}{\Delta}
\newcommand{\e}{\varepsilon}
\newcommand{\g}{\gamma}
\newcommand{\G}{\Gamma}
\renewcommand{\l}{\lambda}
\renewcommand{\L}{\Lambda}

\newcommand{\s}{\sigma}

\renewcommand{\th}{\theta}
\renewcommand{\O}{\Omega}
\renewcommand{\o}{\omega}
\newcommand{\z}{\zeta}

\renewcommand{\i}{\infty}

\newcommand{\bQ}{{\mathbb Q}}
\newcommand{\bE}{{\mathbb E}}

\newcommand{\bP}{{\mathbb P}}
\newcommand{\bR}{{\mathbb R}}

\newcommand{\bZ}{{\mathbb Z}}

\newcommand{\cA}{{\mathcal A}}

\newcommand{\cP}{{\mathcal P}}

\newcommand{\cL}{{\mathcal L}}

\begin{document}

\title{Regular Flows for Diffusions with Rough Drifts}
\author{Fraydoun Rezakhanlou\thanks{This work is supported in part by NSF Grant DMS-1106526.} \\
UC Berkeley \\
Department of Mathematics \\
Berkeley, CA\ \ 94720-3840}

\maketitle

\begin{abstract}

According to DiPerna-Lions theory, velocity fields with weak derivatives in $L^p$ spaces possess weakly regular flows. When a velocity field is perturbed by a white noise, the corresponding
(stochastic) flow is far more regular in spatial variables; a $d$-dimensional
diffusion with a drift in $L^{r,q}$ space ($r$ for the spatial variable and $q$ for the temporal variable)
possesses weak derivatives with stretched exponential bounds, provided that 
$r/d+2/q<1$.  
As an application we show that a Hamiltonian system that is perturbed by a white noise produces a symplectic flow provided that the corresponding Hamiltonian function $H$ satisfies
 $\nabla H\in L^{r,q}$ with $r/d+2/q<1$.  As our second application we derive a Constantin-Iyer type circulation formula for certain weak solutions of Navier-Stokes equation. 
\end{abstract}

\section{Introduction}
\label{sec1}

The velocity field of an incompressible inviscid fluid is modeled by Incompressible Euler Equation
 \begin{equation}\label{eq1.1}
u_t + (Du)u + \nabla P = 0,\ \ \ \nabla \cdot u=0,
\end{equation}
where $u:\bR^d\times [0,\i)\to\bR^d$ represents the velocity field and $P:\bR^d\times [0,\i)\to\bR$
is the pressure. Here and below we write $Du$ and $\nabla P$ for the $x$-derivatives of the vector field $u$ and the scalar-valued function $P$ respectively.
 In the Lagrangian formulation of the fluid, we interpret $u$ as the velocity of generic fluid particles
 and its flow $X:\bR^d\times [0,\i)\to\bR^d$, defined by
\begin{equation}\label{eq1.2}
\frac {d}{dt}X(a,t) = u(X(a,t),t),\ \ \ X(a,0)=a ,
\end{equation}
plays a crucial role in understanding the regularity of solutions of the equation \eqref{eq1.1}.
Since a solution of \eqref{eq1.1} could be singular,
 we need to examine the regularity of the flow $X$ of ordinary differential equations associated with 
rough vector fields. 
Classically, a Lipschitz continuous vector field $u$ results in a Lipschitz flow. 
In a prominent work [DL], DiPerna and Lions constructed a unique flow for \eqref{eq1.2} provided that
$u\in W^{1,p}$ and $\nabla\cdot u\in L^\i$, for some $p\ge 1$. In 2004, Ambrosio [A] extended this result to the case of a vector field $u$ of bounded variation. Recently DeLellis and Crippa [CD] obtained a logarithmic control on the $L^p$--modulus
of continuity of the flow in spatial variable provided that $p>1$.

In the case of an incompressible viscid fluid, the velocity field $u$ satisfies
 the celebrated Navier-Stokes  equation 
\begin{equation}\label{eq1.3}
u_t +(Du)u + \nabla p(x,t) = \nu   \D u, \ \ \ \nabla\cdot u=0.
\end{equation}
 In the corresponding Lagrangian description,  a fluid  particle motion is now modeled by
 a stochastic differential equation (SDE) of the form
\begin{equation}\label{eq1.4}
dX=u(X,t)\ dt+\s dB,
\end{equation}
where $\s=\sqrt{2\nu }$ and $B$ denotes the standard Brownian motion. Since the regularity of solutions to Navier-Stokes equation is a long-standing open problem, 
 we would like to study the regularity of the stochastic flow of
SDE \eqref{eq1.4}  and use such regularity to study \eqref{eq1.3}.
 As it turns out, the flow of SDE \eqref{eq1.4} is far more regular than its 
inviscid analog \eqref{eq1.2}. To state the main result of this article, let us define
\[
\|f\|_{{r,q}}:=\|f\|_{L^{r,q}}:=\left[\int_0^T \left(\int_{\bR^d}|f(x,t)|^r dx\right)^{q/r}dt\right]^{1/q}=
\left[\int_0^T \|f(\cdot,t)\|_{L^r(\bR^d)}^{q}dt\right]^{1/q}.
\]
The space of functions with $\|u\|_{r,q}<\i$ is denoted by $L^{r,q}$. We write $\bP$ and 
$\bE$ for the probability measure and expectation associated with SDE \eqref{eq1.4}. 

\begin{theorem}\label{th1.1} Assume that $\s>0$ and $u\in L^{r,q}$ for some $q\in (2,\i], r\in(d,\i]$, satisfying
\begin{equation}\label{eq1.5}
\d_1:=\frac 12-\frac d{2r}-\frac 1q>0.
\end{equation}
Then SDE \eqref{eq1.4} has a flow $X$ that is weakly differentiable 
with respect to the spatial variable. Moreover, there exist positive constants $C_1=C_1(r,q)$
and $C_0=C_0(r,q)$ such that for every $p\ge 1$,
\begin{equation}\label{eq1.6}
\sup_a\bE\left[\left|D_aX(a,t)\right|^{p}+\left|\left(D_aX(a,t)\right)^{-1}\right|^{p} \right]    \le C_0^p
\exp\left(C_1\s^{-\frac 1{\d_1}\left(1+\frac dr\right)}p^{\frac 1{\d_1}}\|u\|_{r,q}^{\frac 1{\d_1}}
\ t\right).
\end{equation}
\end{theorem}

The following consequence of Theorem~1.1 allows us to go beyond the $p$-th moment and gives an almost Lipschitz regularity of the flow in the spatial variable.

\begin{corollary}\label{cor1.1}
There exist positive constants $C_1'=C_1'(r,q)$ and $C_2=C_2(r,q;\ell)$ such that 
\begin{align}\label{eq1.7}
\sup_a\bP\left(\left|D_aX(a,t)\right|\ge \l\right)\le 
&\exp\left(-C'_1\s^{\frac 1{1-\d_1}\left(1+\frac dr\right)}\|u\|_{r,q}^{-\frac 1{1-\d_1}}
\ t^{-\frac{\d_1}{1-\d_1}}\left(\log^+\frac{\l}{C_0}\right)^{\frac 1{1-\d_1}}\right),\\
\bE \sup_{\substack{|a-b|\le\d\\|a|,|b|\le \ell}}|X(a,t)-X(b,t)|
\le &C_2\ \d \exp\left(C_2\left(\log\frac{\ell}{\d}\right)^{1-\d_1} \s^{-\left(1+\frac dr\right)}\|u\|_{r,q}\ t^{\d_1}\right),
\label{eq1.8}
\end{align}
for every $u$ and $X$ as in Theorem~1.1.
\end{corollary}

As another application of \eqref{eq1.6}, we can show that the flow is jointly 
H\"older continuous in both $x$ and $t$ variables. Define
\[
S_{T,\ell}(X;\d):=\sup_{\substack{|s-t|\le\d\\
0\le s,t\le T}}\sup_{\substack{|a-b|\le\d\\|a|,|b|\le \ell}}{|X(a,t)-X(b,t)|}
\]
\begin{corollary}\label{cor1.2}
For every $\a\in(0,1/2)$, there exists  a constant $C'_2=C_2'(r,q;\ell,T;\a)$ such that 
\begin{equation}\label{eq1.9}
\bE\ S_{T,\ell}(X;\d)
\le C'_2\ \d^\a\exp\left(C'_2 \s^{-\frac 1{\d_1}\left(1+\frac dr\right)}\|u\|_{r,q}^{\frac 1{\d_1}}\right),
\end{equation}
for every $\d>0$, and $u$ and $X$ as in Theorem~1.1.
\end{corollary}

\bigskip
\noindent
{\bf{Remark 1.1}} Clearly our bounds \eqref{eq1.6} and \eqref{eq1.7} are vacuous when $\d_1=0$.
Nonetheless we conjecture that some variants of these bounds would still be true when $\d_1=0$,
 or even when $\d_1<0$. Though we do not expect to have bounds that are uniform in $a$. 
\qed

One of our main motivation behind Theorem 1.1 is its potential applications in Symplectic Topology.
 It also allows us to formulate Navier-Stokes Equation geometrically. To explain this note that
\eqref{eq1.6} allows us to make sense of the pull-back $X_t^*\b$, for any differential form $\b$, where $X_t(\cdot)=X(\cdot,t)$.
  Let us explain this in the case of a $1$-form.
If $\b=w\cdot dx$, or equivalently $\b(x;v)=w(x)\cdot v$,
then
\[
X_t^*\b(x;v)=\b(X_t(a);D_aX_t(a)v)=\left(D_aX_t(a)\right)^*w(X_t(a))\cdot v,
\]
is all well defined. In fact if $w\in L^\i_{loc}$ , then $X^*_tw:=(DX_t)^*\big( w\circ X_t\big)
\in L^p_{loc}$ for every 
$p\in[1,\i)$. In the case that $w\in C^2$, we can make sense
of $\cA_u\b$, where $\cA_u=\cL_u+\nu \D$ with $\cL_u$ denoting the Lie derivative and
\[
\D\left(\sum_i w^i\ dx^i\right)=\sum_i(\D w^i)\ dx^i.
\]
The following theorem explains the role of the operator $\cA_u$.

 \begin{theorem}\label{th1.2} Let $X$ be the flow of SDE \eqref{eq1.4} with $u\in L^{r,q}$ for some $r$ and $q$ satisfying \eqref{eq1.5}. Given $\b^t=w(\cdot,t)\cdot dx$, with $w(\cdot,t)\in C^2$
and $w(x,\cdot)\in C^1$, the process
\begin{equation*}
M^t=X^*_t\b^t-\b^0-\int_0^t X^*_s\left[\dot\b^s+\cA_u\b^s\right]\ ds,
\end{equation*}
is a martingale. (Here $\dot\b^t=w_t(\cdot,t)\cdot dx.$) More precisely,
\begin{equation}\label{eq1.10}
M^t=\int_0^t\sum_{i=1}^d X^*_s\g^s_i\  dB^i(s),
\end{equation}
where $\g^s_i=w_{x^i}(\cdot,s)\cdot dx$.
\end{theorem}

Let us be more precise about the meaning of martingales in our setting. Observe that $M^t$ is a 
$1$-form for each $t$ and we may regard $M^t=M^t(x)$ as a vector-valued function for each $t$.
By Theorem~1.1, this function is locally in $L^p$ for every $p\in[1,\i)$. Now $M^t$ is a martingale in the following sense: If $V(x)$ is a $C^1$ vector field of compact support, then the process
\begin{equation}\label{eq1.11}
 M_t(V):=\int_{\bR^d}M^t(x)\cdot V(x)\ dx,
\end{equation}
is a martingale. As we will see in Section~5, the expression $\cA_u\b^s$ is well-defined weakly; 
only after an integration by parts we can make sense of $M_t(V)$. 
To explain this, recall that by Cartan's formula
\[
\cL_u\b=\hat d i_u\b+i_u \hat d \b,
\]
 where we are writing $\hat d$ for the exterior derivative and $i_u$ denotes the contraction operator in the direction $u$. 
(To avoid confusion with stochastic differential, we use a hat for exterior derivative.) Since $w\in C^1$, we
have no problem to define $i_u \hat d \b$. However we need differentiability of $u$ to make sense
$\hat d i_u\b$ classically. The differentiability of $u$ can be avoided if we integrate against a $C^1$ function
because 
\[
\int \hat d i_u\b(V)\ dx=-\int \b(u)\ (\nabla\cdot V)\ dx.
\]
 
 Let us write
\[
J = \begin{bmatrix}
0 & I_n \\
-I_n & 0
\end{bmatrix} 
\]
where $I_n$ denotes the $n \times n$ identity matrix. 
As a straight forward consequence of Theorem~1.2, we have Corollary 1.1.

\begin{corollary}\label{cor1.2} Assume that $d=2n$, and $u=J\nabla_x H$ in \eqref{eq1.4} for a function
$H:\bR^{d}\x [0,\i)\to\bR$ that is weakly differentiable in $x$-variable and $\nabla_x H\in L^{q,r}$,
for some $q$ and $r$ satisfying  \eqref{eq1.5}. Then the flow $X_t$ is symplectic.
\end{corollary}

\bigskip
\noindent
{\bf{Remark 1.2}} In {\em{Hofer Geometry}} [H], if $\{H_n\}$ is a $L^{\i,1}$-Cauchy sequence of Hamiltonian functions
of compact support,  the the corresponding flows $\{\phi^{H_n}\}$ is 
a Cauchy sequence of symplectic flows with respect to the
{\em{Hoper metric}}. Completion of the group of such symplectic transformations with respect to the Hofer metric is not understood. In view of Corollary~1.3, we may wonder whether or not some kind of a limit exists for the family of the flows $\{X=X^\s:\s>0\}$ as $\s\to 0$.
\qed

As our next application, let us assume that $u$ is a solution of the backward Navier-Stokes Equation:
\begin{equation}\label{eq1.12}
u_t +(Du)u + \nabla P(x,t) + \nu   \D u=0, \ \ \ \nabla\cdot u=0.
\end{equation}
(We use backward equation \eqref{eq1.12} instead of the forward equation \eqref{eq1.3} to simplify our presentation.) A more geometric formulation of \eqref{eq1.12} is achieved by 
writing an equation for the evolution of the $1$-form $\a^t=u(\cdot,t)\cdot dx$:
 \begin{equation}\label{eq1.13}
\dot \a^t+\cA_u\a^t-\hat dL^t=0,
\end{equation}
where $L^t(x)=\frac 12|u(x,t)|^2-P(x,t)$. A natural way to approximate Navier-Stokes Equation is via {\em{Camassa-Holm}}-type equations of the form
\begin{equation}\label{eq1.14}
v_t +(Dv)w + (Dw)^*v-\nabla_x \bar L(x,t) + \nu   \D v=0, \ \ \ \nabla\cdot v=0,\ \ \ w=v\ast_x \z,
\end{equation}
where $\zeta(x)$ is a smooth function. In the classical Camassa-Holm Equation, $v=w-\e \D w$ which leads to
$u=v\ast_x \z^\e$. In this case both $v=v^\e$ and $w=w^\e$ depend on $\e$ and according to a classical result of Foias et al. [FHT],  the sequences $(w^\e,v^\e)$ are precompact in low $\e$ 
limit and if $(u,u)$ is any limit point, then $u$ is a weak-solution of \eqref{eq1.12}. We say $u$ is a $(r,q)$-regular solution of \eqref{eq1.12} if it can be approximated by a sequence of solutions $(v^\e,w^\e)$ of Camassa-Holm equation such that
\begin{equation}\label{eq1.15}
\sup_{\e>0}\|w^\e\|_{r,q}<\i.
\end{equation}

 \begin{theorem}\label{th1.3}
Let $u$ be a $(r,q)$-regular solution of Navier-Stokes Equation \eqref{eq1.12} 
for some $r$ and $q$ satisfying \eqref{eq1.5}. Then for any smooth divergence free vector field $Z$ of compact support, the process
\begin{equation}\label{eq1.16}
\int_{\bR^d}X^*_t\a^t(Z)\ dx,
\end{equation}
is a martingale. Moreover, if $Du\in L^2,$ then the process $X_t^*\hat d\a^t$ is also a martingale.
\end{theorem}

\medskip
\noindent
{\bf{Remark 1.3}}
 According to a classical result of Serrin [S], a weak solution of \eqref{eq1.3} 
is smooth if $u(\cdot,0)\in L^2$ and $u\in L^{r,q}$ for some $r$ and $q$ satisfying \eqref{eq1.5}.
We may also use Theorem~1.3  to show that any $(r,q)$-regular solution is smooth.
If we have equality in \eqref{eq1.5} and $r<\i$, the regularity of solutions can be found in the work of Fabes, Jones and Riviere [FJR]. 
Based on this, it is natural to ask what type of regularity for the flow $X$ is available in the extreme case $\d_1=0$ (see Remark~1.1 above). We leave this for future investigation.
\qed

Here is a short review of various classical and recent results on SDE \eqref{eq1.4}:

\medskip\noindent
{\bf{1.}} Classical Ito's theorem guarantees that \eqref{eq1.4} has a unique (strong) solution if $u$ is Lipschitz
continuous in spatial variable, uniformly in time.

\medskip\noindent
{\bf{2.}} By a yet another classical work of Bismut, Elworthy and Kunita
(see for example [RW] or [K]), \eqref{eq1.4} has a smooth flow 
with smooth inverse if $u$ is smooth.

\medskip\noindent
{\bf{3.}} Zvonkin [Z] in 1974 showed that \eqref{eq1.4} has a unique solution if $d=1$ and $u\in L^{\i,\i}$.
This result was extended to higher dimension by Veretynikov [V] in 1979.

\medskip\noindent
{\bf{4.}} Flandoli et al. [FGP] (2010) have shown that if $u$ is H\"older-continuous of H\"older exponent $\a$
in spatial variable, then the flow $X$ is also H\"older-continuous of H\"older exponent $\a'$ in spatial variable, for any $\a'<\a$.

\medskip\noindent
{\bf{5.}} Fedrizzi and Flandoli [FF] (2010) establish $X\in W^{1,p}$ for every $p\ge 2$,
provided that $u\in L^{r,q}$ for some $r$ and $q$ satisfying \eqref{eq1.5}.
 Though no bound on $D_aX$ is given in [FF].

 \medskip\noindent
{\bf{6.}} Mohammad et al. [MNP] (2014) establish $\bE|D_aX(a,t)|^p<\i$ for every $p\ge 1$
 provided that $u\in L^{\i,\i}$.
\qed

An important ingredient for the work of Mohammad et al. is a bound of Davie [D]
 (see Theorem~2.1) that works
for $u\in L^{\i,\i}$. In this paper we adopt [MNP] approach and achieve Theorem 1.1 by generalizing
Davie's bound to the case $u\in L^{r,q}$ with $r$ and $q$ satisfying \eqref{eq1.5}. 
In fact Davie proves such a bound by reducing it to a certain double integral. 
It is worth mentioning that such a reduction is applicable only if we assume a stronger condition 
\begin{equation}\label{eq1.17}
\d_2:=\frac 14-\frac d{2r}-\frac 1q>0,
\end{equation}
We refer to Subsection~4.2 for more details.

The organization of the paper is as follows: 
\begin{itemize}
\item In Section~2 we establish Theorem~1.1 and its corollaries, assuming that a Davie-type bound 
(Theorem~2.1) is available under the assumption $\d_1>0$.
\item In Section~3 we reduce the proof of Theorem~2.1 to bounding certain block-type integrals
(Theorem 3.1). 
\item Section~4 is devoted to the proof of Theorem~3.1.
\item In Section~5 we discuss symplectic diffusions and prove Theorems~1.2 and 1.3.
\end{itemize}

\section{Proof of Theorem 1.1 and Its Corollaries}
\label{sec2}
 \setcounter{equation}{0}

As a preparation for the proof of Theorem~1.1, we state one theorem and two lemmas.
We write $x^1,\dots,x^d$ for coordinates of $x$ and $f_{x^i}$ 
for the partial derivative of $f$ with respect to $x^i$.
\begin{theorem}\label{th2.1} For every $r$ and $q$ satisfying \eqref{eq1.5},
there exists a constant $C_3=C_3(r,q)$
 such that for any continuously differentiable functions 
$b^1,\dots,b^n:\bR^d\times [0,1]\to \bR$ of compact support, 
and indices $\a_1,\dots,\a_n\in\{1,\dots,d\}$, we have
\begin{equation}\label{eq2.1}
\left|\bE\int_{\D^n}\prod_{i=1}^{n}b^i_{x^{\a_i}}(a+\s B(t_i),t_i)dt_i\right|
\le C_3^n\s^{-n\left(\frac dr+1\right)} n^{-n\d_1}
t^{n\d_1}\prod_{i=1}^{n}\|b^i\|_{q,r},
\end{equation}
where
\[
\Delta^n=\D^n(t)=\{(t_1,\dots,t_n)\ :\ 0\le t_1\le \dots \le t_n\le t\}.
\]
\end{theorem}

\begin{lemma}\label{lem2.1}
For every $r$ and $q$ satisfying \eqref{eq1.5},
there exists a constant $C_4=C_4(r,q)$
  such that
\begin{equation}
\label{eq2.2}\sup_a\bE \ \exp\left[\l\int_0^{t}|u|^2(a+\s B(s),s)\ ds\right]\le C_4 
\exp\left[C_4\s^{-\frac{d}{r\d_1}}\l^{\frac 1{2\d_1}}
\|u\|_{r,q}^{\frac 1{\d_1}}\ t\right].
\end{equation}
\end{lemma}

\begin{lemma}\label{lem2.2} For every $\b\in(0,1)$ and $p>(d+1)\b^{-1}$, we can find
a constant $C_5=C_5(p,\b)$
  such that for every continuous function $X\in\bR^d\times[0,T]\to\bR^d$,
\begin{equation}
\label{eq2.3}
S_{T,\ell}(X;\d)\le C_5\d^{\b-\frac{d+1}p}\ \left\{\int_{s,t\in[0,T]}\int_{|x|,|y|\le \ell}
\frac{|X(x,t)-X(y,t)|^p}{|(x,t)-(y,s)|^{\b p+d+1}}\  dxdydtds\right\}^{\frac 1p}.
\end{equation}
\end{lemma}

Lemma~2.2 is the known as Garsia-Rodemich-Rumsey Inequality and its proof can be found in
[SV]. Inequality \eqref{eq2.2} is a Khasminskii type bound and its proof will be given at the end of this section. Theorem~2.1 is the main ingredient for the proof of Theorem~1.1
 and was established by Davie when $q=r=\i$. The proof of Theorem~2.1 will be given in the next section.

Before embarking on the proof of Theorem~1.1, let us outline our strategy. 
\begin{itemize}
\item{\bf{(i)}} We first assume that $u$ is a smooth function of compact support. This guarantees
 that the flow $X$ is a diffeomorphism in spatial variables. For such a drift $u,$ we establish \eqref{eq1.6}. Note that the right-hand side of
\eqref{eq1.6} depends only $\|u\|_{r,q}$ norm and is independent of the smoothness of $u$.
\item{\bf{(ii)}} Given $u\in L^{r,q}$ with $(r,q)$ satisfying \eqref{eq1.5},
we choose a sequence of smooth functions $\{u_N\}_{N=1}^{\i}$ of compact supports such that $\|u_N-u\|_{r,q}\to 0$ in large $N$ limit. Writing $\cP_N$ for the law of the corresponding flow $X^N$, we use 
Corollary~1.2 to show that the sequence 
$\{\cP_N\}_{N=1}^{\i}$ is tight. Then by standard arguments we can show that any limit point of 
$\{\cP_N\}_{N=1}^{\i}$ is a law of a flow $X$ that satisfies \eqref{eq1.4} and that
the bounds \eqref{eq1.6}-\eqref{eq1.9} are valid.
\end{itemize}

 For the proof of \eqref{eq1.6} we follow 
[MNP] closely; Step~1 and part of Step~2 are almost identical to the proof of Lemma~7 in [MNP]. 

\bigskip
\noindent
{\bf{Proof of Theorem~1.1. Step 1.}} We prove \eqref{eq1.6} assuming that $u$ is a smooth vector field of compact support as in Part (i) of the above outline. We leave Part (ii) for Section~5 where Theorem~1.2 is established. From
\[
dX(a,t)=u(X(a,t),t) \ dt+\s dB,\ \ \ \ X(a,t)=a,
\]
we can readily deduce 
\begin{equation}\label{eq2.4}
\frac {d}{dt}D_aX(a,t)=D_aX_t(a,t)=D_xu(X(a,t),t)D_aX(a,t),\ \ \ \ D_aX(a,0)=I,
\end{equation}
where $I$ denotes the $d\x d$ identity matrix. Regarding \eqref{eq2.4} as an ODE for
$D_aX(a,t)$, this equation has a unique solution and this solution is given by
\begin{equation}\label{eq2.5}
D_aX(a,t)=I+\sum_{n=1}^\i \int_{\Delta^n(t)}
D_xu(X(a,t_n),t_n)\dots D_xu(X(a,t_1),t_1)\ d{t_1}\dots d{t_n},
\end{equation}
provided that this series is convergent. ($\D^n$ was defined right after \eqref{eq2.1}.) 

As for the inverse $(D_aX)^{-1}$, observe 
\[
\frac {d}{dt}\left(D_aX(a,t)\right)^{-1}=-\left(D_aX(a,t)\right)^{-1}\frac {d}{dt}\left(D_aX(a,t)\right)
\left(D_aX(a,t)\right)^{-1}.
\]
This and \eqref{eq2.4} yields
\begin{equation*}
\frac {d}{dt}\left(D_aX(a,t)\right)^{-1}=-\left(D_aX(a,t)\right)^{-1}D_xu(X(a,t),t),\ \ \ \
 \left(D_aX(a,0)\right)^{-1}=I.
\end{equation*}
Regarding this as an ODE for
$(D_aX(a,t))^{-1}$, this equation has a unique solution and this solution is given by
\begin{equation}\label{eq2.6}
D_aX(a,t)=I+\sum_{n=1}^\i (-1)^{-n}\int_{\Delta^n(t)}D_xu(X(a,t_1),t_1)
\dots D_xu(X(a,t_n),t_n)\ d{t_1}\dots d{t_n},
\end{equation}
provided that this series is convergent.

We use \eqref{eq2.5} to bound $\big|D_aX(a,t)\big|$. (In the same fashion, we may use 
\eqref{eq2.6} to bound $\big|(D_aX(a,t))^{-1}\big|$.)
This is achieved by bounding the summand in
\eqref{eq2.5}, which in the end verifies the convergence of the series and the validity of \eqref{eq2.5}.

Using the matrix norm $|[a_{ij}]|=\sum_{i,j}|a_{ij}|$, we have
\begin{equation}\label{eq2.7}
\left[\bE|D_aX_t(a,t)|^p\right]^{\frac 1p}\le d+\sum_{n=1}^\i A_n^{\frac 1p},
\end{equation}
where
\[
A_n=\bE\left| \int_{\D^n}
D_xu(X(a,t_n),t_n)\dots D_xu(X(a,t_1),t_1)\ d{t_1}\dots d{t_n}\right|^p.
\] 
Writing $x=(x^1,\dots,x^d)$ and $u=(u^1,\dots,u^d)$, we may assert
\begin{equation}\label{eq2.8}
A_n\le d^{(n+1)(p-1)}\sum_{i_0,\dots,i_{n}=1}^dA_n(i_0,\dots,i_{n}),
\end{equation}
where $A_n(i_0,\dots,i_{n})$ is given by
\begin{equation*}
\bE\left| \int_{\D^n}
u^{i_0}_{x^{i_1}}(X(a,t_n),t_n)u^{i_1}_{x^{i_2}}(X(a,t_{n-1}),t_{n-1})\dots u^{i_{n-1}}_{x^{i_{n}}}(X(a,t_1),t_1)\ d{t_1}\dots d{t_n}\right|^p
\end{equation*}
On the other hand, for $p$ an even integer, we can drop absolute values and
 express $A_n(i_0,\dots,i_{n})$
as a sum of at most $p^{np}$ terms of the form $B_{np}(j_1,k_1,\dots,
j_{np},k_{np}),$ that is given by
\begin{equation*}
\bE\int_{\D^{np}}
u^{k_1}_{x^{j_1}}(X(a,s_{np}),s_{np})\dots
 u^{k_{np}}_{x^{j_{np}}}(X(a,s_{1}),s_{1})\ d{s_1}\dots d{s_{np}},
\end{equation*}
for $j_1,k_1,\dots,j_{np},k_{np}\in\{1,\dots.d\}$. This is because there are at most $p^{np}$ many ways
to form $s_1\le \dots\le s_{np}$ out of $p$ many groups of the form
\[
\{ t^i_1\le\dots\le t^i_n\},\ \ \ \ i=1,\dots,p.
\]
(Once $s_1\le \dots\le s_{\ell}$ are selected, there are at most $p$ many possibilities for our next selection
$s_{\ell+1}$.)

{\bf{Step 2.}} Writing $\bQ^a$ for the law of $(a+\s B(s):s\in[0,t])$ with $B(\cdot)$ representing  a standard Brownian motion
that starts from $0$, and applying Girsanov's formula, we may write $B_{np}$ as
\[
\int\left[\int_{\Delta^{np}}
u^{k_1}_{x^{j_1}}(x(s_{np}),s_{np})\dots
 u^{k_{np}}_{x^{j_{np}}}(x(s_{1}),s_{1})\ d{s_1}\dots d{s_{np}}\right]\ 
M(x(\cdot))\ \bQ^a(dx(\cdot)),
\]
where
\[
M(x(\cdot))=M_u(x(\cdot))=\exp\left(\frac 1{2\nu}\int_0^t u(x(s),s)\cdot dx(s)-
\frac1{4\nu}\int_0^t |u(x(s),s)|^2\ ds\right).
\]
This, by Schwartz' inequality, is bounded above by 
\begin{equation*}
D_{np}(j_1,k_1,\dots,
j_{np},k_{np})^{\frac 12}\ \left(\int M^2 \ d\bQ^a\right)^{\frac 12},
\end{equation*}
where $D_{np}(j_1,k_1,\dots,j_{np},k_{np})$ is given by
\[
\int\left[\int_{\D^{np}}
u^{k_1}_{x^{j_1}}(x(s_{np}),s_{np})\dots
 u^{k_{np}}_{x^{j_{np}}}(x(s_{1}),s_{1})\ d{s_1}\dots d{s_{np}}\right]^2\ \bQ^a(dx(\cdot)).
\]
As in Step 1, we may express $D_{np}(j_1,k_1,\dots,j_{np},k_{np})$  as a sum of at most $2^{np}$
many terms of the form $E_{2np}(j_1,k_1,\dots,j_{2np},k_{2np})$, that are defined as  
\[
\int\left[\int_{\D^{2np}}
u^{k_1}_{x^{j_1}}(x(s_{2np}),s_{2np})\dots
 u^{k_{2np}}_{x^{j_{2np}}}(x(s_{1}),s_{1})\ d{s_1}\dots d{s_{2np}}\right]\ \bQ^a(dx(\cdot)).
\] 
By Theorem~2.1, 
\begin{equation*}
\big|E_{2np}(j_1,k_1,\dots,j_{2np},k_{2np})\big|\le
C_3^{2np}(2\nu)^{-np\left(\frac dr+1\right)}(2np)^{-2np\d_1}t^{2np\d_1}\|u\|^{2np}_{r,q}.
 \end{equation*}
This in turn implies
\begin{equation}\label{eq2.9}
\big|D_{np}(j_1,k_1,\dots,j_{np},k_{np})\big|\le 2^{np}
C_3^{2np}(2\nu)^{-np\left(\frac dr+1\right)}(2np)^{-2np\d_1}t^{2np\d_1}\|u\|^{2np}_{r,q}.
 \end{equation}
Furthermore, by Lemma~2.1,
\begin{align}\nonumber
\int M_u^2 \ d\bQ^a&=\int \exp\left(\frac 1{\nu}\int_0^t u(x(s),s)\cdot dx(s)-
 \frac 1{2\nu}\int_0^t |u(x(s),s)|^2\ ds\right) \ d\bQ^a\\ \nonumber
&=\int \left(M_{4u}\right)^{\frac 12}\exp\left(\frac 3{2\nu }\int_0^t |u(x(s),s)|^2\ ds\right) \ d\bQ^a\\
&\le\left(\int M_{4u}\ d\bQ^a\right)^{\frac 12}
\left(\int \exp\left(\frac 3\nu \int_0^t |u(x(s),s)|^2\ ds\right) \ d\bQ^a\right)^{\frac 12}
\label{eq2.10}\\
&=\left(\int \exp\left(\frac 3\nu \int_0^t |u(x(s),s)|^2\ ds\right) \ d\bQ^a\right)^{\frac 12}\nonumber\\
&\le C_4\exp\left(c_0\ \nu^{-\frac 1{2\d_1}\left(1+\frac dr\right)}\|u\|_{r,q}^{\frac 1{\d_1}}
\ t\right),
\nonumber
\end{align}
where for the third equality we use the fact that $M_{4u}$ is a $\bQ^a$-martingale. From \eqref{eq2.7}--  \eqref{eq2.10} we deduce that for positive even integer $p$,
\begin{equation}\label{eq2.11}
\left[\bE|D_aX(a,t)|^p\right]^{\frac 1p}\le d+C_4 Z(p)
\exp\left(c_02^{-1}p^{-1}\nu^{-\frac 1{2\d_1}\left(1+\frac dr\right)}\|u\|_{r,q}^{\frac 1{\d_1}}
\ t\right)
\end{equation}
where
\begin{align*}
Z(p)=&\sum_{n=1}^\i 
d^{n+1}2^{\frac n2}p^nC_3^{n}\nu^{-\frac n2\left(\frac dr+1\right)}(2np)^{-n\d_1}t^{n\d_1}\|u\|^{n}_{r,q}\\
&\le \sum_{n=1}^\i  
\left(c_1\nu^{-\frac 1{2\d_1}\left(\frac dr +1\right)}p^{\frac 1{\d_1}-1}\ t\|u\|_{r,q}^{\frac 1{\d_1}}\right)^{n\d_1}
n^{-n\d_1}\\ 
&=:\sum_{n=1}^\i W^{n\d_1}n^{-n\d_1}.
\end{align*}
On the other hand, 
\begin{align*}
\sum_{n=1}^\i W^{n\d_1}n^{-n\d_1}&=\sum_{n=1}^\i \left((2W)^n n^{-n}\right)^{\d_1}2^{-n\d_1}
\le c_2 \left(\sum_{n=1}^\i(2W)^n n^{-n}\ 2^{-n\d_1}\right)^{\d_1}\\
&\le c_3 \left(\sum_{n=1}^\i(2W)^n(n!)^{-1}\ 2^{-n\d_1}\right)^{\d_1}\le c_3e^W,
\end{align*}
where we used Stirling's formula for the second inequality. 
From this and \eqref{eq2.11} we deduce,
\begin{equation}\label{eq2.12}
\bE|D_aX(a,t)|^p\le 
\exp\left(c_4p+c_4\nu^{-\frac 1{2\d_1}\left(1+\frac dr\right)}p^{\frac 1{\d_1}}\|u\|_{r,q}^{\frac 1{\d_1}}
\ t\right).
\end{equation}
The bound \eqref{eq2.12} is true for every even integer $p\ge 2$. By changing the constant $c_4$
if necessary, we can guarantee that it is also true for every real $p\in[1,\i)$. As we mentioned earlier, with a verbatim argument we can establish the analog of \eqref{eq2.12} for $(D_aX)^{-1}$.
\qed

\bigskip
\noindent
{\bf{Proof of Corollary 1.1.}} We start with the proof of \eqref{eq1.7}. From \eqref{eq1.6} and
Chebyshev's inequality we learn that for every $p,\l\in(1,\i)$,
\begin{equation*}
\bP\left(|D_aX(a,t)|\ge \l\right)\le \l^{-p}C_0^pe^{Ap^{\frac 1{\d_1}}}=\exp\left(-p\log\l+p\log C_0
+Ap^{\frac 1{\d_1}}\right),
\end{equation*}
where  
\begin{equation}\label{eq2.13}
		A=C_1\s^{-\frac 1{\d_1}\left(1+\frac dr\right)}\|u\|_{r,q}^{\frac 1{\d_1}}\ t.
\end{equation}
We optimize this bound by choosing $\log\l=A\d_1^{-1}p^{\frac 1{\d_1}-1}+\log C_0$;
\begin{equation*}
\bP\left(|D_aX(a,t)|\ge \l\right)\le C_0e^{A\left(1-\frac 1{\d_1}\right)p^{\frac 1{\d_1}}}
\le  \exp\left(-C'_1A^{-\frac {\d_1}{1-\d_1}}\left(\log^+\frac{\l}{C_0}\right)^{\frac 1{1-\d_1}}\right),
\end{equation*}
for a positive constant $C'_1$. This completes the proof of \eqref{eq1.7}.

We next turn to the proof of \eqref{eq1.8}. Set
\[
\o_\ell(\d)=\sup_{\substack{|a-b|\le\d\\|a|,|b|\le \ell}}|X(a,t)-X(b,t)|.
\]
By Morrey's inequality [E],
\begin{equation*}
\o_\ell(\d)\le c_0\d^{1-\frac dp}\left(\int_{|z|\le 2\ell}|D_aX(z,t)|^p\ dz\right)^{\frac 1p},
\end{equation*}
for every $p>d$, and for a universal constant $c_0$ that can be chosen to be independent of $p$.
This, \eqref{eq1.6} and H\"older's inequality imply
\[
\bE \o_\ell(\d)\le c_0\d^{1-\frac dp}\left(\int_{|z|\le 2\ell}\bE|D_aX(z,t)|^p\ dz\right)^{\frac 1p}
\le c_1\d^{1-\frac dp}\ell^{\frac dp}e^{Ap^{\frac 1{\d_1}-1}},
\]
with $A$ defined by \eqref{eq2.13}. We optimize this bound by choosing
\[
\log\frac {\ell^d}{\d^d}=A\left(\frac 1{\d_1}-1\right)p^{\frac 1{\d_1}}.
\]
For such a choice of $p$ we deduce 
\begin{equation*}
\bE \o_\ell(\d)
\le c_1\d \exp\left(c_2\left(\log\frac{\ell}{\d}\right)^{1-\d_1} A^{\d_1}\right),
\end{equation*}
for a positive constant $c_2$. This completes the proof of \eqref{eq1.8}.
\qed

\bigskip
\noindent
{\bf{Proof of Corollary~1.2.}} We use \eqref{eq1.6} to assert that for $t\in[0,T]$,
\begin{align}\nonumber
\bE|X(x,t)-X(y,t)|^p&= \bE \left|\int_0^1 D_aX(\th x+(1-\th)y,t)\cdot(x-y)\ d\th\right|^p\\
&\le  C_0^p\exp\left(C_1\nu^{-\frac 1{2\d_1}\left(1+\frac dr\right)}p^{\frac 1{\d_1}}\|u\|_{r,q}^{\frac 1{\d_1}}T\right)|x-y|^{p}.\label{eq2.14}
\end{align}
On the other hand, by Girsanov's formula and
H\"older's inequality,
\begin{align*}
\bE|X(y,t)-X(y,s)|^p&=  \int{|x(t)-x(s)|^p}\ M(x(\cdot))\ \bQ^y(dx(\cdot))\\
&\le  \left(\int{|x(t)-x(s)|^{p\g}}\ \bQ^y(dx(\cdot))\right)^{\frac 1{\g}} \left(\int M^{\g'}\ d\bQ^y\right)^{\frac 1{\g'}}\\
&\le  c_0|t-s|^{\frac p2}\  \left(\int M^{\g'}\ d\bQ^y\right)^{\frac 1{\g'}}\\
&\le c_1\exp\left(c_2\ \nu^{-\frac 1{2\d_1}\left(1+\frac dr\right)}\|u\|_{r,q}^{\frac 1{\d_1}}
\ T\right)\ |t-s|^{\frac p2},
\end{align*}
where $1/\g+1/\g'=1$, $\bQ^y$ and $M$ were defined in the beginning of Step~2 of the proof of Theorem~1.1, and for the last inequality we follow \eqref{eq2.10}. 
From this and \eqref{eq2.14} we deduce
\[
\bE|X(x,t)-X(y,s)|^p\le c_3^p\exp\left(c_3\ \nu^{-\frac 1{2\d_1}\left(1+\frac dr\right)}
\|u\|_{r,q}^{\frac 1{\d_1}}\ T\right)\left(|x-y|^p+|t-s|^{\frac p2}\right).
\]
This and Lemma~2.2 imply,
\begin{equation*}
\bE\ S_{T,\ell}(X;\d)^p\le c_4c_3^p\exp\left(c_3\ \nu^{-\frac 1{2\d_1}\left(1+\frac dr\right)}
\|u\|_{r,q}^{\frac 1{\d_1}}\ T\right)\ \d^{\b p-d-1},
\end{equation*}
with $c_4<\i$ if $\b\in(0,1/2)$. (Here we have used $|x-y|^{p}\le c|x-y|^{p/2}$ for $|x|,|y|\le \ell$.)
Finally we choose $p$ so that $\b-(d+1)/p=\a$ to complete the proof.
\qed

\bigskip
We end this section with the proof of Lemma~2.1. Let us make some preparations.
We write $p(x,t)=(t\nu)^{-d/2}p(x/ \sqrt{t\nu})$ with 
\begin{equation}\label{eq2.15}
p(z)=(4\pi )^{-d/2}\exp\left(-|z|^2/4\right).
\end{equation}
Throughout the paper we need to bound $L^r$ norms of $p(\cdot,s)$ and its spatial derivatives. These bounds are stated in Lemma 2.3 below. The elementary proof of this lemma is omitted.
\begin{lemma}\label{lem2.3}
For every $r\in[1,\i]$ and nonnegative integer $k$, there exists a constant $C_5(k,r)$ such that if $\tilde p(\cdot,s)$ denotes a $k$-th spatial derivative of 
$p(\cdot,s)$, then 
\begin{equation}\label{eq2.16}
\|\tilde p(\cdot,s)\|_{L^{r'}}\le C_5(s\nu)^{-\frac d{2r}-\frac k2},
\end{equation}
 for every $s>0$,where $r'=r/(r-1)$.
\end{lemma}
 We are now ready to establish \eqref{eq2.2}.

\bigskip
\noindent
{\bf{Proof of Lemma~2.1.}} The proof is based on Khasminskii's trick. We first show
that there exists a constant $c_1$ such that
\begin{equation}\label{eq2.17}
\sup_a\sup_{\th\in[0,t]}\bE \int_0^{t}|u|^2(a+\s B(s),s+\th)\ ds\le c_1 \nu^{-d/r}t^{2\d_1}
\|u\|_{r,q}^2.
\end{equation}
This is a straight forward consequences of H\"older's inequality:
\begin{align*}
\int_0^t\int |u|^2 (a+x,s+\th)p(x,s)\ dx ds&\le 
\int_0^t\left(\int |u|^r(x,s)\ dx\right)^{\frac 2r}\left(\int p^{r'}(x,s+\th)
\ dx\right)^{\frac 1{r'}}\ ds\\
&\le c_0\int_0^t\left(\int |u|^r(x,s)\ dx\right)^{\frac 2r}(s\nu)^{-\frac d2 (1-\frac 1{r'})}\ ds\\
&= c_0 \|u\|_{r,q}^2\left(\int_0^t (s\nu)^{-\frac d2 (1-\frac 1{r'})q'}\ ds\right)^{\frac 1{q'}}\\
&= c_1 \nu^{-d/r}t^{2\d_1}\|u\|_{r,q}^2,
\end{align*}
where $\frac 2r+\frac 1{r'}=1$ and $\frac 2q+\frac 1{q'}=1$. Given $\l>0$, choose $t_0$ such that 
\[
c_1 \nu^{-d/r}t_0^{2\d_1}\|u\|_{r,q}^2\l=\frac 12=:\a_0,
\]
and use Khasminskii's trick (see for example [S]) to deduce 
\[
\sup_a\bE \ \exp\left[\l\int_0^{t_0}|u|^2(a+B(s),s)\ ds\right]\le (1-\a_0)^{-1}=2, 
\]
from \eqref{eq2.17}. This and Markov property yields
\[
\sup_a\bE \ \exp\left[\l\int_0^{\ell t_0}|u|^2(a+B(s),s)\ ds\right]\le 2^\ell. 
\]
This implies \eqref{eq2.2} after choosing $\ell=[t/t_0]+1$.
\qed

\section{Proof of Theorem~2.1}
\label{sec3}
 \setcounter{equation}{0}

The main ingredient for the proof of Theorem~2.1 is a bound on certain {\em{block integrals}}. 
In this section we state a crucial bound for block integrals and show how such bounds can be used to establish Theorem~2.1.

We say a function $h$ is of type $j$ if it is  a spatial partial $j$th-derivative of $p$.
We can readily show that if $h$ is of type $j$, then
\begin{equation}\label{eq3.1}
h(z,s)\le C_6 \ (\nu s)^{-\frac j2}p(z,2s),
\end{equation}
for a constant $C_6$. Define
\[
 \D^k=\D^k(t_0,t)=\{(t_1,\dots,t_k):\ t_0\le t_1\le \dots\le t_k\le t\}.
\]
For our purposes, we would like to bound {\em{block integrals}} $I^k(f_1,\dots,f_{k})$, 
where 
\begin{itemize}
\item $I^1(f_1)=\int_{\D^1}\int_{\bR^{d}} f_1(z_1,t_1)p^{(1)}(z_1,t_1-t_0)(t-t_1)^\a\ dz_1dt_1,$ with 
$p^{(1)}$ being of type $1$.
\item $I^2(f_1,f_2)$ is defined as 
\[
\int_{\D^2}\int_{\bR^{2d}} f_1(z_1,t_1)f_2(z_2,t_2)p(z_1,t_1-t_0)
 p^{(2)}(z_2-z_1,t_2-t_1) (t-t_2)^\a\ dz_1dz_2dt_1dt_2,
\]
with 
$p^{(2)}$ being of type $2$.
\item For  $k>2$, we define $I^k(f_1,\dots,f_{k})$ by
 \begin{align*}
\int_{\D^k}\int_{\bR^{kd}} f_1(z_1,t_1)&p(z_1,t_1-t_0)\left\{\prod_{i=2}^{k-1}f_i(z_i,t_i)
p^{(1)}_i(z_i-z_{i-1},t_i-t_{i-1})\ dz_idt_i\right\}\\
&\ \ \ \times  f_k(z_k,t_k)p^{(2)}(z_{k}-z_{k-1},t_k-t_{k-1}) (t-t_k)^\a\ dz_1dz_kdt_1dt_k,
\end{align*}
 where $p^{(2)}$ is of type $2$ and $p^{(1)}_i$ is of type $1$ for 
$i=2,\dots,k-1$. 
\end{itemize}
Our main result on block integrals is Theorem~3.1.
\begin{theorem}\label{thm3.1} There exists a constant $C_7=C_7(r,q)$ such that
\begin{equation}\label{eq3.2}
|I_k(f_1,\dots,f_{k})|\le C_7^k\nu^{-k\left(\frac d{2r}+\frac 12\right)}\g_k(\a)(t-t_0)^{\a+k\d_1}
\prod_{i=1}^k\|f_i\|_{r,q} ,
\end{equation}
where
\[
\g_k(\a)=\frac{\a^\a }{(\a+k\d_1)^{\a+k\d_1}}.
\]
(By convention, $0^0=1$.)
\end{theorem}

Armed with Theorem~3.1, we are now ready to give a proof for \eqref{eq2.1}

 \bigskip
\noindent
{\bf{Proof of Theorem~2.1.}} {\bf{Step 1.}} Let us write $R$ for the left-hand side of \eqref{eq2.1}.
We certainly have
\begin{align*}
R&=\left|
\int_{\D^n}\int_{\bR^{dn}}\prod_{i=1}^{n}b^i_{x^{\a_i}}(a+y_i,t_i)
p(y_i-y_{i-1},t_i-t_{i-1})dy_idt_i\right|\\
& =\left|\int_{\D^n}\int_{\bR^{dn}}
\prod_{i=1}^{n}b^i_{y_i^{\a_i}}(a+y_i,t_i)
p(y_i-y_{i-1},t_i-t_{i-1})dy_idt_i\right|,
\end{align*}
where  $y_0=0$, $y_i^{\a_i}$ denotes the $\a_i$-th coordinate of $y_i\in\bR^d$, and $p$ was defined by
\eqref{eq2.15}.
After some integration by parts we learn 
\begin{equation}\label{eq3.3}
R=\left|\sum_{r=1}^{2^{n-1}}\e_rI(\b_1(r),\dots,\b_n(r))\right|,
\end{equation}
where each $\e_r$ is either $1$ or $-1$, the indices $\b_1(r),\dots,\b_n(r)$ are in $\{0,1,2\}$ and satisfy
$\sum_i\b_i(r)=n$,
and the expression $I$ has the form
\[
I(\b_1,\dots,\b_n)=\int_{\D^n}\int_{\bR^{dn}}\prod_{i=1}^{n}b^i(a+y_i,t_i)
q^{\b_i}(y_i-y_{i-1},t_i-t_{i-1})dy_idt_i.
\]
 Here $q^0(a,t)=p(a,t)$, $q^1(a,t)=p_{a^j}(a,t)$ for some 
$j\in\{1,\dots d\}$ and $q^2(a,t)=p_{a^ja^k}(a,t)$ for some 
$j,k\in\{1,\dots d\}$. Recall that by $p_{a^j}$ and $p_{a^ja^k}$, we mean partial derivatives with respect to coordinates $a^j$ and $a^j,a^k$ respectively. As a result, \eqref{eq2.1} would follow, if we can find a constant $c_1$ such that for all $\b_1,\dots,\b_n$,
\begin{equation}\label{eq3.4}
\left|I(\b_1,\dots,\b_n)\right|
\le c_1^n  \nu^{-\frac n2\left(\frac dr+1\right)} \left(n\d_1\right)^{-n\d_1}t^{n\d_1}
\prod_{i=1}^{n}\|b^i\|_{q,r}.
\end{equation}
By induction on $n$, we can readily show that the type of $n$-tuple
$(\b_1,\dots,\b_n)$ that appears in \eqref{eq3.3} can be decomposed into blocks of sizes $n_1,\dots,n_\ell$ 
such that if
\[
m_0=0, \ \ m_1=n_1,\ \ m_2=n_1+n_2,\ \dots,\  m_\ell=n_1+n_2+\dots+n_\ell=n,
\]
then each block $(\b_{m_{i-1}+1},\dots,\b_{m_i})$ satisfies the following conditions: 
\begin{itemize}
\item
If $n_i=1$, then $\b_{m_{i-1}+1}=1$.
\item
If $n_i=2$, then $\b_{m_{i-1}+1}=0$ and $\b_{m_{i-1}+2}=\b_{m_i}=2$.
\item
If $n_i>2$, then $\b_{m_{i-1}+1}=0$ and $\b_{m_i}=2$ and all $\b_s$ in between are $1$.
\end{itemize}

{\bf{Step 2.}} When $\ell>1$, we set 
 \[
J_{\ell-1}(t_{m_{\ell-1}},y_{m_{\ell-1}})=\int_{\D^{n_{\ell}}(t_{m_{\ell-1}},t)}
\int_{\bR^{dn_{\ell}}}\prod_{i=m_{\ell-1}+1}^{n}b^i(a+y_i,t_i)
q^{\b_i}(y_i-y_{i-1},t_i-t_{i-1})\ dy_idt_i.
\]
In the case of $\ell>2$, we inductively define
 \[
J_{j}(t_{m_j},y_{m_j})=\int_{\D^{n_{j+1}}(t_{m_j},t)}\int_{\bR^{dn_{j+1}}}J_{j+1}(t_{m_{j+1}},y_{m_{j+1}})\prod_{i=m_j+1}^{m_{j+1}}b^i(a+y_i,t_i)
q^{\b_i}(y_i-y_{i-1},t_i-t_{i-1})\ dy_idt_i.
\]
for $j=\ell-2,\dots,1.$ This allows us to write
\[
I(\b_1,\dots,\b_n)=\int_{\D^{n_1}(0,t)}\int_{\bR^{dn_1}}J_{1}(t_{m_{1}},y_{m_{1}})\prod_{i=1}^{m_{1}}b^i(a+y_i,t_i)
q^{\b_i}(y_i-y_{i-1},t_i-t_{i-1})\ dy_idt_i.
\]
We then apply Theorem~3.1 to assert  
\[
\left|J_{\ell-1}(t_{m_{\ell-1}},y_{m_{\ell-1}})\right|
\le C^{n_{\ell}}_7 \nu^{-n_{\ell}\left(\frac d{2r}+\frac 12\right)}
\left(n_{\ell}\d_1\right)^{-n_{\ell}\d_1}\left(t-t_{m_{\ell-1}}\right)^{n_{\ell}\d_1}
\prod_{i=m_{\ell-1}+1}^{n}\|b^i\|_{r,q}.
\]
This allows us to express 
\[
J_{\ell-1}(t_{m_{\ell-1}},y_{m_{\ell-1}})=\hat J_{\ell-1}(t_{m_{\ell-1}},y_{m_{\ell-1}})\left(t-t_{m_{\ell-1}}\right)^{n_{\ell}\d_1},
\]
with $\hat J_{\ell-1}$ satisfying
\[
\left|\hat J_{\ell-1}(t_{m_{\ell-1}},y_{m_{\ell-1}})\right|
\le C^{n_{\ell}}_7 \nu^{-n_{\ell}\left(\frac d{2r}+\frac 12\right)}
\left(n_{\ell}\d_1\right)^{-n_{\ell}\d_1}
\prod_{i=m_{\ell-1}+1}^{n}\|b^i\|_{r,q}.
\]
After replacing $b^{m_{\ell-1}}$ with $b^{m_{\ell-1}}\hat J_{\ell-1}$, we apply Theorem~3.1 again 
to assert
\begin{align*}
|J_{\ell-2}(t_{m_{\ell-2}},y_{m_{\ell-2}})|\le C^{n_{\ell}+n_{\ell-1}}_7 
&\nu^{-\left(n_{\ell}+n_{\ell-1}\right)\left(\frac d{2r}+\frac 12\right)}
\left(\left(n_{\ell}+n_{\ell-1}\right)\d_1\right)^{-\left(n_{\ell}+n_{\ell-1}\right)\d_1}\\
&\ \  \x\left(t-t_{m_{\ell-2}}\right)^{(n_{\ell}+n_{\ell-1})\d_1}\prod_{i=m_{\ell-2}+1}^{n}\|b^i\|_{r,q},
\end{align*}
provided that $\ell>2$. Continuing this inductively we arrive at \eqref{eq3.4} for $c_1=C_7$.
The bound \eqref{eq3.4} in turn implies \eqref{eq2.1} for $C_3=2C_7\d_1^{-\d_1}$.
\qed

\section{Bounding Block Integrals}
\label{sec4}
\setcounter{equation}{0}

\subsection{Proof of Theorem~3.1}

As preparation for the proof of Theorem~3.1, we  establish two lemmas.
 The first lemma is a slight generalization of \eqref{eq3.2} when $k=2$. 
Given $\b\ge 0$, define $I'(f_1,f_2)$ by
\[
\int_{\D^2(t)}\int_{\bR^{2d}} f_1(z_1,t_1)f_2(z_2,t_2)p(z_1,2t_1)
 p^{(2)}(z_2-z_1,t_2-t_1)t_1^\b (t-t_2)^\a\ dz_1dz_2dt_1dt_2,
\]
where $\D^2=\D^2(0,t)$. Also define $J(f_1,\dots,f_\ell;z_{\ell+1},t_{\ell+1})$ by
\[
\int_{\D^\ell(t_{\ell+1})}\int_{\bR^{d\ell}}p(z_1,t_1)f_1(z_1,t_1)
p_{\ell+1}^{(1)}(z_{\ell+1}-z_{\ell},t_{\ell+1}-t_{\ell})
\prod_{i=2}^{\ell }f_i(z_i,t_i)
p^{(1)}_i(z_i-z_{i-1},t_i-t_{i-1})\ \prod_{i=1}^{\ell}dz_idt_i ,
\]
where
$$\D^\ell(t_{\ell+1})=\{(t_1,\dots,t_{\ell}):0\le t_1\le \dots\le t_\ell\le t_{\ell+1}\}.$$

\begin{lemma}\label{lem4.1} There exists a constant $C_8=C_8(r,q)$ such that for $\a,\b\ge 0$,
\begin{equation}\label{eq4.1}
|I'(f_1,f_{2})|\le C_8\nu^{-\left(\frac d{r}+1\right)}\zeta(\a,\b)t^{2\d_1+\a+\b}
\|f_1\|_{r,q}\|f_2\|_{r,q} .
\end{equation}
where
\[
\zeta(\a,\b)=\frac{\a^{\a} \b^{\b}(\b+1)^{\d_1-\frac d{2r}}}{(\a+\b+2\d_1)^{\a+\b+2\d_1}}.
\]
 (By convention $0^0=1$.)
\end{lemma}

\begin{lemma}\label{lem4.2}There exists a constant $C_9=C_9(r,q)$ such that
\begin{equation}\label{eq4.2}
|J(f_1,\dots,f_\ell;z_{\ell+1},t_{\ell+1})|\le C_9^\ell \
\nu^{-\left(1+\frac dr\right) \frac \ell 2}\ (\ell\d_1)^{-\ell\d_1}
\ p\left(z_{\ell+1},2t_{\ell+1}\right)t_{\ell+1}^{\ell \d_1}\prod_{i=1}^\ell\|f_i\|_{r,q}. 
\end{equation} 
\end{lemma}

Before embarking on the proofs of Lemmas~4.1 and 4.2, let us 
recall the relationship between generalized Beta function and Gamma function $\G$.
\begin{lemma}\label{lem4.3} For every $\a_0\dots\a_n>0$,
\begin{equation}\label{eq4.3}
\int_{t_0\le t_1\le \dots\le t_{n+1}}\prod_{i=0}^{n}(t_{i+1}-t_i)^{\a_i-1}\ \prod_{i=1}^{n}dt_i=
(t_{n+1}-t_0)^{\sum_{i=0}^n\a_i-1}\ \frac{\prod_{i=0}^{n}\G(\a_i)}{\G\left(\sum_{i=0}^n\a_i\right)}.
\end{equation}
\end{lemma}

The elementary proof of Lemma~4.3 is omitted.

\bigskip\noindent
{\bf{Proof of Lemma 4.1.}} 
 We write $\D$ for the set
$\D^2=\D^2(t)$ and set 
\[
f_i'(s)=\|f_i(\cdot,s)\|_{L^r(\bR^d)}.
\]
{\bf {Step 1.}}
We decompose 
$I'=I_1+I_2$ where $I_i$ is obtained from $I$ by replacing the domain of integration $\D=\D^2$ with 
$\D_i.$ The sets $\D_1$ and $\D_2$ are defined by
\begin{align*}
\D_1&=\{(t_1,t_2)\in \D:\  t_1\le t_2-t_1\},\\
\D_2&=\{(t_1,t_2)\in \D:\  t_2-t_1\le t_1\}.
\end{align*}
The term $I_1$ is easily bounded with the aid of our $L^r$ bounds on $p$ and $p^{(2)}$: If
we set
\[
\eta(\a,\b;q):=\left(\frac {\G(\b q'+\d_1q')\G(\d_1q')\G(\a q' +1)}{\G(2\d_1q'+(\a+\b) q'+1)}\right)^{\frac 1{q'}},
\]
with $q'=q/(q-1)$, then by Lemma~2.3 and
H\"older's inequality, the expression $|I_1|$ is bounded above by
\begin{align*}
 c_0&\nu^{-\left(\frac dr+1\right)}\int_{\D_1}f_1'(t_1)f_2'(t_2)
 t_1^{\b-\frac d{2r}}(t_2-t_1)^{-\frac d{2r}-1}(t-t_2)^{\a }\ dt_1dt_2\\
&\le c_0\nu^{-\left(\frac dr+1\right)}\|f_1\|_{r,q}\|f_2\|_{r,q}\ 
\left(\int_{\D_1} t_1^{\left(\b-\frac {d}{2r}\right)q'}(t_2-t_1)^{-\left(\frac d{2r}+1\right)q'}(t-t_2)^{\a q'}
\ dt_1dt_2\right)^{\frac 1{q'}}\\
&\le c_0\nu^{-\left(\frac dr+1\right)}\|f_1\|_{r,q}\|f_2\|_{r,q}\ 
\left(\int_{\D_1} t_1^{\left(\b-\frac d{2r}-\frac 12\right)q'}
(t_2-t_1)^{-\left(\frac d{2r}+\frac 12\right)q'}(t-t_2)^{\a q'}\ dt_1dt_2\right)^{\frac 1{q'}}\\
&= c_0\nu^{-\left(\frac dr+1\right)}\|f_1\|_{r,q}\|f_2\|_{r,q}\ t^{2\d_1+\a+\b}\eta(\a,\b;q),
\end{align*}
 where for the second inequality
we used the fact that $t_2-t_1\ge t_1$ in the set $\D_1$, and for the equality we used
the fact 
\[
\left(\frac d{2r}+\frac 12\right)q'<1,
\]
which is the same as \eqref{eq1.5}. In summary,
\begin{equation}\label{eq4.4}
|I_1|\le c_0\nu^{-\left(\frac dr+1\right)}\eta(\a,\b;q)\|f_1\|_{r,q}\|f_2\|_{r,q}\ t^{2\d_1+\a+\b},
\end{equation}
It remains to bound $I_2$.

{\bf{Step 2.}} We next decompose $I_2$ as $I_{21}+I_{22}$, where $I_{21}$ is obtained from
$I_2$ by restricting the domain of $dz_1dz_2$-integration to a set of points $(z_1,z_2)$ such that
 $|z_2-z_1|/\sqrt{\nu t_1}$ stays away from zero.
Though this restriction is done so that the product structure of 
$f_1(z_1,t_1)f_2(z_2,t_2)$ is not destroyed. For this 
purpose, we decompose $\bR^{2d}$ into cells 
$$B_{k\ell}(t_1,t_2)=B_k(t_1)\x B_\ell(t_2),$$ 
where for $k=(k^1,\dots,k^d)$, the set
$B_k(s)$ denotes the set of $z=(z^1,\dots,z^d)$ such that
\[
z^i/\sqrt{\nu s}\in[k^i,k^i+1),
\]
for $i=1,\dots,d$. We now write $|k-\ell|_1=\sum_{i=1}^d|k^i-\ell^i|$ for the $L^1$ distance between
$k,\ell\in \bR^d$, and set 
\[
I_{21}=\sum_{(k,\ell)\in\L_1}I_{2}(k,\ell),\ \ \ I_{22}=\sum_{(k,\ell)\in\L_2}I_{2}(k,\ell),
\]
where
\begin{align*}
\L_1&=\left\{(k,\ell):\ k,\ell\in \bZ^d,\
 |k|_1\notin\left[|\ell|_1-4d,\sqrt 2|\ell|_1+4d\right]\right\}\\
\L_{2}&=\left\{(k,\ell):\ k,\ell\in \bZ^d,\
 |k|_1\in\left[|\ell|_1-4d,\sqrt 2|\ell|_1+4d\right]\right\},
\end{align*}
and $I_{2}(k,\ell)$ is defined by
\[
\int_{\D_2}\int_{B_{k\ell}(t_1,t_2)} f_1(z_1,t_1)f_2(z_2,t_2)p(z_1,2t_1)
 p^{(2)}(z_2-z_1,t_2-t_1)t_1^\b(t-t_2)^\a\ dz_1dz_2dt_1dt_2.
\]
To bound $I_{21}$, assume that $(z_1,z_2)\in B_{k\ell}$ for some $k,\ell$ satisfying
either $|k|_1> \sqrt 2|\ell|_1+4d$, or $|\ell|_1> |k|_1+4d$. If the former occurs and $(t_1,t_2)\in \D_2$, 
then 
\begin{align*}
|z_2-z_1|_1&=|(z_2-\ell\sqrt{\nu t_2})-(z_1-k\sqrt{\nu t_1})+\ell\sqrt{\nu t_2}-k\sqrt{\nu t_1}|_1\\
&\ge |\ell\sqrt{\nu t_2}-k\sqrt{\nu t_1}|_1-d(\sqrt{\nu t_2}+\sqrt{\nu t_1})\\
&\ge|k|_1\sqrt{\nu t_1}- |\ell|_1\sqrt{\nu t_2}-d(\sqrt{\nu t_2}+\sqrt{\nu t_1})\\
&\ge \left(\sqrt 2|\ell|_1+4d\right)\sqrt{\nu t_1}-|\ell|_1\sqrt{\nu t_2}-d(\sqrt{\nu t_2}+\sqrt{\nu t_1})\\
&\ge \left(\sqrt 2|\ell|_1+4d\right)\sqrt{\nu t_1}-|\ell|_1\sqrt{2\nu t_1}-d(\sqrt{2\nu t_1}+\sqrt{\nu t_1})\\
&\ge d \sqrt{\nu t_1}.
\end{align*}
If the latter occurs and $(t_1,t_2)\in \D_2$, 
then 
\begin{align*}
|z_2-z_1|_1&=|(z_2-\ell\sqrt{\nu t_2})-(z_1-k\sqrt{\nu t_1})+\ell\sqrt{\nu t_2}-k\sqrt{\nu t_1}|_1\\
&\ge |\ell|_1\sqrt{\nu t_2}-|k|_1\sqrt{\nu t_1}-d(\sqrt{\nu t_2}+\sqrt{\nu t_1})\\
&\ge \left(|k|_1+4d\right)\sqrt{\nu t_1}-|k|_1\sqrt{\nu t_1}-d(\sqrt{\nu t_2}+\sqrt{\nu t_1})\\
&\ge 3d\sqrt{\nu t_1}-d\sqrt{\nu t_2}\ge d\sqrt{\nu t_1}.
\end{align*}
In any case, we always have
\[
|z_2-z_1|^2\ge d^{-1} |z_2-z_1|_1^2\ge d{\nu t_1}.
\]
From this and \eqref{eq3.1} we learn
\begin{align*}
| p^{(2)}(z_2-z_1,t_2-t_1)|&\le \frac {c_1}{\nu(t_2-t_1)}p(z_2-z_1,2(t_2-t_1))\\
&= \frac {c_1}{\sqrt{\nu(t_2-t_1)}} 
\frac {|z_2-z_1|}{\sqrt{\nu(t_2-t_1)}}\frac 1{|z_2-z_1|}p(z_2-z_1,2(t_2-t_1))\\
&\le \frac {c_2}{\sqrt{\nu(t_2-t_1)}} \frac 1{|z_2-z_1|}p(z_2-z_1,4(t_2-t_1))\\
&\le \frac {c_2}{\sqrt{\nu(t_2-t_1)}} \frac 1{\sqrt{d\nu t_1}}p(z_2-z_1,4(t_2-t_1)).
\end{align*}
This and Lemma~2.3 imply that the term $I_{21}$ is bounded above by a constant multiple of 
\begin{align*}
\nu^{-1}&\int_{\D}\int_{\bR^{2d}}| f_1(z_1,t_1)f_2(z_2,t_2)|t_1^{\b-\frac 12}p(z_1,2t_1)|
 (t_2-t_1)^{-\frac 12}p(z_2-z_1,4(t_2-t_1))\\
&\ \ \ \ \ \ \ \ \ \ \ \ \ \ \ \ \ \ \ \ \ \ \ \ \times (t-t_2)^\a\ dz_1dz_2dt_1dt_2\\
&\le c_3\nu^{-\left(\frac dr+1\right)}\int_{\D}f_1'(t_1)f_2'(t_2) t_1^{\b-\frac d{2r}-\frac 12}
(t_2-t_1)^{-\frac d{2r}-\frac 12}(t-t_2)^{\a }\ dt_1dt_2\\
&\le c_3\nu^{-\left(\frac dr+1\right)}\|f_1\|_{r,q}\|f_2\|_{r,q}\ 
\left(\int_{\D} t_1^{\left(\b-\frac {d}{2r}-\frac 12\right)q'}(t_2-t_1)^{-(\frac d{2r}+\frac 12)q'}(t-t_2)^{\a q'}
\ dt_1dt_2\right)^{\frac 1{q'}}\\
&=c_3\nu^{-\left(\frac dr+1\right)}\|f_1\|_{r,q}\|f_2\|_{r,q}\ t^{2\d_1+\a+\b}\eta(\a,\b;q).
\end{align*}
In summary,
\begin{equation}\label{eq4.5}
|I_{21}|\le c_{4}\nu^{-\left(\frac dr+1\right)}\eta(\a,\b;q)
\|f_1\|_{r,q}\|f_2\|_{r,q}\ t^{2\d_1+\a+\b},
\end{equation}
It remains to bound $I_{22}$.

{\bf{Step 3.}} Let us write
\begin{align*}
f_{1k}(z_1,t_1)&=f_1(z_1,t_1)1\!\!1(z_1\in B_k(t_1))p(z_1,2t_1),\\
f_{2\ell}(z_2,t_2)&=f_2(z_2,t_2)1\!\!1(z_2\in B_\ell(t_2)),
\end{align*}
so that $I_2(k,\ell)$ can be expressed as
\[
\int_{\D_2}\int_{\bR^{2d}} f_{1k}(z_1,t_1)f_{2\ell}(z_2,t_2)p^{(2)}(z_2-z_1,t_2-t_1)t_1^\b(t-t_2)^\a\ dz_1dz_2dt_1dt_2.
\]
Recall that $ p^{(2)}$ is a function of type $2$. This means that  $p^{(2)}(z,s)=p_{z^iz^j}(z,t)$ is a second derivative of $p$. By Plancheral's formula we 
learn that  $I_2(k,\ell)$ equals to
\[
-(2\pi)^2\int_{\D_2}\int_{\bR^{d}}\xi^i\xi^j \hat f_{1k}(\xi,t_1)\check f_{2\ell}(\xi,t_2)
e^{-4\pi^2\nu(t_2-t_1)|\xi|^2}t_1^\b(t-t_2)^\a\ d\xi  dt_1dt_2,
\]
where 
\[
\hat h(\xi,s)=\int e^{-2i\pi x\cdot\xi}h(x,s)\ dx,\ \ \  \check h(\xi,s)=\int e^{2i\pi x\cdot\xi}h(x,s)\ dx.
\]
As a result,  the term $|I_2(k,\ell)|$ is bounded above by
\begin{align*}
 2 {\pi^2}\int_{\D_2}\int_{\bR^{d}}& \left( \d^{-1}(\nu t_1)^{\frac d2}|\hat f_{1k}(\xi,t_1)|^2+
\d (\nu t_1)^{-\frac d2}|\check f_{2\ell}(\xi,t_2)|^2\right)\ |\xi|^2\\
&\ \ \ \ \ \ \ \ \ \times e^{-4\pi^2\nu(t_2-t_1)|\xi|^2}t_1^\b(t-t_2)^\a\ d\xi  dt_1dt_2
=:I_2^1(k,\ell)+I_2^2(k,\ell).
\end{align*}
for any $\d>0$. Further, $I_2^1(k,\ell)$ equals to
\begin{align*}
 \frac {2\pi^2}{\d}&\int_{\D_2}\int_{\bR^{d}}(\nu t_1)^{\frac d2}|\hat f_{1k}(\xi,t_1)|^2
 |\xi|^2 e^{-4\pi^2\nu(t_2-t_1)|\xi|^2}t_1^\b(t-t_2)^\a\ d\xi dt_1dt_2\\
&\le \frac {2\pi^2}{\d}\int_{\D}\int_{\bR^{d}}(\nu t_1)^{\frac d2}|\hat f_{1k}(\xi,t_1)|^2
 |\xi|^2 e^{-4\pi^2\nu(t_2-t_1)|\xi|^2}t_1^\b(t-t_1)^\a\ d\xi dt_1dt_2\\
&\le \frac {1}{2\nu\d}\int_{0}^t\int_{\bR^{d}}(\nu t_1)^{\frac d2}|\hat f_{1k}(\xi,t_1)|^2t_1^\b(t-t_1)^\a
\ d\xi dt_1\\
&=\frac {1}{2\nu\d}\int_{0}^t\int_{\bR^{d}}(\nu t_1)^{\frac d2}| f_{1k}(x,t_1)|^2t_1^\b(t-t_1)^\a
\ dx dt_1\\
&\le  {c_{5}e^{-\frac 14(|k|-\sqrt d)^{+2}}}(\nu \d )^{-1}\int_{0}^t 
(\nu t_1)^{-\frac d2}\int_{B_k(t_1)}| f(x,t_1)|^2t_1^\b(t-t_1)^\a
\ dx dt_1,
\end{align*}
where we used $\a\ge 0$ and $t_1\le t_2$ for the first inequality.
We now apply H\"older's inequality to assert that $I_2^1(k,\ell)$ is bounded above by
\begin{align*}
c_{6}&e^{-\frac 14(|k|-\sqrt d)^{+2}}{(\nu \d) ^{-1}}\int_{0}^t 
\left((\nu t_1)^{-\frac d2}\int_{B_k(t_1)}| f(x,t_1)|^r\ dx\right)^{\frac 2r}t_1^\b(t-t_1)^\a
\  dt_1\\
&\le{c_{7}e^{-\frac 15|k|^2}}{(\nu \d) ^{-1}}\int_{0}^t f'(t_1)^2\
(\nu t_1)^{-\frac dr}\ t_1^\b(t-t_1)^\a
\  dt_1\\
&\le{c_{7}\nu^{-\left(\frac dr+1\right)}e^{-\frac 15|k|^2}}{\d ^{-1}}\|f_1\|^2_{r,q}\left(\int_{0}^t 
t_1^{-\frac {dq}{r(q-2)}}t_1^{\frac {\b q}{q-2}}(t-t_1)^{\frac {\a q}{q-2}}
\  dt_1\right)^{\frac {q-2}q}\\
&={c_{7}\nu^{-\left(\frac dr+1\right)}e^{-\frac 15|k|^2}}{\d ^{-1}}\|f_1\|^2_{r,q}\ \eta'(\a,\b)
\ t^{2\d_1+\a+\b},
\end{align*}
where
 \[
\eta'(\a,\b;q)=\left(
\frac{\G((2\d_1+\b)q'')\G(\a q''+1)}{\G((2\d_1+\a+\b)q''+1)}\right)^{\frac 1{q''}},
\]
for $q''=q/(q-2)$. Note that for the equality, we have used the fact that $dq/\big(r(q-2)\big)<1$,
which is equivalent to \eqref{eq1.5}.
In summary,
\begin{equation}\label{eq4.6}
|I_{2}^1(k,\ell)|\le c_{7}\nu^{-\left(\frac dr+1\right)}
e^{-\frac 15|k|^2}{\d ^{-1}}\eta'(\a,\b;q)\|f_1\|_{r,q}^2\ t^{2\d_1+\a+\b}.
\end{equation}
On the other hand, $I_2^2(k,\ell)$ is bounded above by
\begin{align*}
 2  \pi^2\d&\int_{\D_2}\int_{\bR^{d}}(\nu t_1)^{-\frac d2}|\check f_{2\ell}(\xi,t_2)|^2
|\xi|^2 e^{-4\pi^2\nu(t_2-t_1)|\xi|^2}t_1^\b(t-t_2)^\a\ d\xi dt_1dt_2\\
&\le 2 {\pi^2\d}2^{d/2}\int_{\D_2}\int_{\bR^{d}}(\nu t_2)^{-\frac d2}|\check f_{2\ell}(\xi,t_2)|^2
|\xi|^2 e^{-4\pi^2\nu(t_2-t_1)|\xi|^2}t_2^\b(t-t_2)^\a\ d\xi dt_1dt_2\\
&\le   \d\nu^{-1}2^{d/2-1}\int_{0}^t\int_{\bR^{d}}(\nu t_2)^{-\frac d2}|\check f_{2\ell}(\xi,t_2)|^2\
t_2^\b(t-t_2)^\a
\ d\xi dt_2\\
&=  \d \nu^{-1}2^{d/2-1}
\int_{0}^t\int_{\bR^{d}}(\nu t_2)^{-\frac d2}| f_{2\ell}(x,t_2)|^2\ t_2^\b(t-t_2)^\a
\ dx dt_2\\
&=  \d\nu^{-1}2^{d/2-1}\int_{0}^t\int_{B_\ell(t_2)}(\nu t_2)^{-\frac d2}| f_2(x,t_2)|^2t_2^\b(t-t_2)^\a
\ dx dt_2\\
&\le c_{8} \d\nu^{-1}\int_{0}^t (\nu t_2)^{-\frac {d}{r}}f_2'(t_2)^2t_2^\b(t-t_2)^\a
 dt_2\\
&\le  c_{8}\d\nu^{-\left(\frac dr+1\right)}
\|f_2\|_{r,q}^2 \left(\int_{0}^t t_2^{\left(\b-\frac {d}{r}\right)q''}(t-t_2)^{\a q''}
 dt_2\right)^{\frac {1}{q''}}\\
&=  c_{8}\d\nu^{-\left(\frac dr+1\right)}\|f_2\|_{r,q}^2 \eta'(\a,\b;q)t^{2\d_1+\a+\b},
\end{align*}
where we used $t_1\le t_2\le 2t_1$ for the first inequality.
 In summary,
\begin{equation}\label{eq4.7}
|I_{2}^2(k,\ell)|\le c_{8}{\d }\nu^{-\left(\frac dr+1\right)}\eta'(\a,\b;q)\|f_2\|_{r,q}^2\ t^{2\d_1+\a+\b}.
\end{equation}
We choose $\d=e^{-|k|^2/10}$ and use \eqref{eq4.6} and \eqref{eq4.7} to deduce
\begin{equation*}
|I_2(k,\ell)|=|I_{1}^2(k,\ell)+I_{2}^2(k,\ell)|
\le c_{9}e^{-\frac{|k|^2}{10}}\nu^{-\left(\frac dr+1\right)}\eta'(\a,\b;q)
\left[\|f_1\|_{r,q}^2+\|f_2\|_{r,q}^2\right]\ t^{2\d_1+\a+\b}.
\end{equation*}
From this and the definition of $I_{22}$ we learn
\begin{align}\nonumber
|I_{22}|&\le \sum_{(k,\ell)\in\L_2}|I_{2}(k,\ell)|\\ \label{eq4.8}
&\le c_{10}\nu^{-\left(\frac dr+1\right)}\eta'(\a,\b;q)
\left[\|f_1\|_{r,q}^2+\|f_2\|_{r,q}^2\right]\ t^{2\d_1+\a+\b}\sum_{k}|k|^de^{-\frac{|k|^2}{10}}\\
&= c_{11}\nu^{-\left(\frac dr+1\right)}\eta'(\a,\b;q)
\left[\|f_1\|_{r,q}^2+\|f_2\|_{r,q}^2\right]\ t^{2\d_1+\a+\b}.\nonumber
\end{align}

{\bf{Final Step.}} From \eqref{eq4.4}, \eqref{eq4.5} and \eqref{eq4.8} we learn that there exists 
a constant $c_{12}$ such that if $\|f_1\|_{r,q},\|f_2\|_{r,q}\le 1$, then
\[
| I'(f_1,f_2)|\le c_{12}\nu^{-\left(\frac dr+1\right)}\left[\eta+\eta'\right](\a,\b;q)\ t^{2\d_1+\a+\b}.
\]
From this and a scaling argument we deduce that for every $f_1$ and $f_2$,
\begin{equation}\label{eq4.9}
| I'(f_1,f_2)|\le c_{12}\nu^{-\left(\frac dr+1\right)}\left[\eta+\eta'\right](\a,\b;q)\ t^{2\d_1+\a+\b}
\|f_1\|_{r,q}\|f_2\|_{r,q}.
\end{equation}
First assume that $\b\ge 1$. 
Using Stirling's formula, we can readily show that there exist constants $c_{13}$ and $c_{14}$ such that
\begin{align*}
\eta(\a,\b;q)&\le c_{13}  \frac{\a^{\a+\frac1{2q'}} (\b+\d_1-1+q^{-1})^{\b+\d_1-\frac1{2q'}}}{(\a+\b+2\d_1)^{\a+\b+2\d_1+\frac1{2q'}}}\le c_{14}\frac {\a^{\a}\b^{\b}}
{( \a+\b+2\d_1)^{\a+\b+2\d_1}},\\
\eta'(\a,\b;q)&\le c_{13} \frac{\a^{\a+\frac1{2q''} }(\b+2\d_1-1+2q^{-1})^{\b+2\d_1-\frac1{2q''}}}{(\a+\b+2\d_1)^{\a+\b+2\d_1+\frac1{2q''}}}\le c_{14}\zeta(\a,\b),
\end{align*}
because 
\[
\d_1-(2q')^{-1}, \d_1-1-q^{-1},2\d_1-1+2q^{-1}<0,
\]
 This and \eqref{eq4.9} imply \eqref{eq4.1} when $\b\ge 1$. The case $\b\in[0,1)$ can be treated likewise.
\qed

\bigskip
\noindent
{\bf{Proof of Lemma 4.2.}} As before, we define $r'$ and $q'$ by $r'=r/(r-1)$ and $q'=q/(q-1)$.
By H\"older's inequality we have 
\[
|J(f_1,\dots,f_\ell;z_{\ell+1},t_{\ell+1})|\le \prod_{i=1}^\ell\|f_i\|_{r,q} 
\left(\int_{\D^\ell(t_{\ell+1})} A(t_1,\dots ,t_{\ell};z_{\ell+1},t_{\ell+1})^{\frac {q'}{r'}}
\prod_{i=1}^\ell dt_i\right)^{\frac 1{q'}},
\]
where $A(t_1,\dots ,t_{\ell};z_{\ell+1},t_{\ell+1})$ equals
\begin{equation*}
 \int_{\bR^{d\ell}}p(z_1,t_1)^{r'}\prod_{i=2}^{\ell +1}
\left|p^{(1)}_i(z_i-z_{i-1},t_i-t_{i-1})\right|^{r'}\ \prod_{i=1}^\ell dz_i .
\end{equation*}
On the other hand, by completing squares (or Markov Property) we know
\begin{equation*}
\int_{\bR^{d}}p(z,s)p(a-z,t)\ dz=p(a,s+t)\int_{\bR^{d}}p\left(z-\frac t{s+t}a,\frac {st}{s+t}\right)\ dz
=p(a,s+t).
\end{equation*}
From this and \eqref{eq3.1} we deduce that $A(t_1,\dots ,t_{\ell};z_{\ell+1},t_{\ell+1})$
is bounded above by 
\begin{align*}
 &C_6^{\ell-1}\prod_{i=2}^{\ell+1 }(\nu(t_{i}-t_{i-1}))^{-\frac{r'}2}  \int_{\bR^{d\ell}} p(z_1,t_1)^{r'}\prod_{i=2}^{\ell +1}
p(z_i-z_{i-1},2(t_i-t_{i-1}))^{r'}\ \prod_{i=1}^\ell dz_i\\
&\ \ \le 2^{\frac {dr'}2} C_6^{\ell-1}\prod_{i=2}^{\ell+1 }(\nu(t_{i}-t_{i-1}))^{-\frac{r'}2} \int_{\bR^{d\ell}} p(z_1,2t_1)^{r'}\prod_{i=2}^{\ell +1}
p(z_i-z_{i-1},2(t_i-t_{i-1}))^{r'}\ \prod_{i=1}^\ell dz_i\\
&\ \ \le c_1^{\ell}(\nu  t_1)^{-\frac{dr'}2+\frac d2}\prod_{i=1}^\ell
(\nu(t_i-t_{i-1}))^{-r'\left(\frac d2+\frac 12\right)+\frac d2}\\
&\ \ \ \ \ \ \ \ \ \ \ \ \ \ \ \x \int_{\bR^{d\ell}} p(z_1,2t_1/r')\prod_{i=2}^{\ell +1}
p(z_i-z_{i-1},2(t_i-t_{i-1})/r')\ \prod_{i=1}^\ell dz_i\\
&\ \ = c_1^{\ell}p\left(z_{\ell+1},\frac {2t_{\ell+1}}{r'}\right)
(\nu  t_1)^{-\frac{dr'}2+\frac {dr'}2}\prod_{i=1}^\ell
(\nu(t_{i+1}-t_i))^{-r'\left(\frac d2+\frac 12\right)+\frac d2}.
\end{align*}
Hence $|J(f_1,\dots,f_\ell;z_{\ell+1},t_{\ell+1})|$ is bounded above by
\begin{align*}
c_1^\ell \ \nu^{-\left(1+\frac dr\right)\frac \ell 2}&\  t_{\ell+1}^{\frac d{2r}}
 \ p\left(z_{\ell+1}, 2t_{\ell+1}\right)\ \prod_{i=1}^\ell\|f_i\|_{r,q}\\ 
&\x \left(\int_{\D^\ell(t_{\ell+1})} \ t_1^{-\frac{dq'}2+\frac {dq'}{2r'}}\prod_{i=1}^\ell
(t_{i+1}-t_i)^{-q'\left(\frac d2+\frac 12\right)+\frac {dq'}{2r'}}\
 dt_i\right)^{\frac 1{q'}}\\
=  & c_1^\ell \ \nu^{-\left(1+\frac dr\right)\frac \ell 2}\ t_{\ell+1}^{\frac d{2r}}\ p\left(z_{\ell+1},2t_{\ell+1}\right)
t_{\ell+1}^{-\frac d{2r}+\ell \d_1}\eta_\ell
\ \prod_{i=1}^\ell\|f_i\|_{r,q}\\
=  & c_1^\ell\ \nu^{-\left(1+\frac dr\right)\frac \ell 2}\  p\left(z_{\ell+1},2t_{\ell+1}\right)
t_{\ell+1}^{\ell \d_1}\ \eta_\ell
\ \prod_{i=1}^\ell\|f_i\|_{r,q},
\end{align*}
where
\[
\eta_\ell=\left(\frac{\G(\d_1 q')^\ell\G((\d_1+1/2) q')}
{\G\left(\left(1/2+(\ell+1))\d_1\right) q'\right)}\right)^{\frac 1{q'}}\le c_2^\ell (\ell\d_1)^{-\ell\d_1},
\]
by Stirling's formula. This completes the proof of \eqref{eq4.2}.

\qed

\bigskip\noindent
{\bf{Proof of Theorem~3.1 when $k=1$.}} Without loss of generality, we may assume that $t_0=0$.
In this case, the poof of \eqref{eq3.2} is an immediate consequence of Lemma~2.3, H\"older's inequality and
\eqref{eq4.3}:
\begin{align*}
|I_1(f_1)|&\le c_0\int_0^t f_1'(t_1)(\nu t_1)^{-\left(\frac d{2r}+\frac 12\right)}
(t-t_1)^\a\ dt_1\\
&\le c_0\nu^{-\left(\frac d{2r}+\frac 12\right)}\|f_1\|_{r,q}
\left(\int_{0}^t t_1^{-\left(\frac d{2r}+\frac 12\right)q'}(t-t_1)^{\a q'}
\ dt_1\right)^{\frac 1{q'}}\\
&=c_1\nu^{-\left(\frac d{2r}+\frac 12\right)}\eta(\a)\|f_1\|_{r,q}\ t^{\a+\d_1},
\end{align*}
where
\[
\eta(\a):=\left(\frac {\G(\d_1q')\G(\a q' +1)}{\G(\d_1q'+\a q'+1)}\right)^{\frac 1{q'}},
\]
Using Stirling's formula, we can readily show that $\eta(\a)\le c_2\g_1(\a)$,
for a constant $c_2$. This completes the proof when
$k=1$.
\qed

\bigskip\noindent
{\bf{Proof of Theorem~3.1 when $k\ge 2$.}} Without loss of generality, we may assume that $t_0=0$.
Evidently \eqref{eq4.1} in the case of $\b=0$ implies \eqref{eq3.2} when $k=2$ because $p(z_1,t_1)\le 2^{\frac d2}
p(z_1,2t_1)$ (or we can readily show that Lemma~4.1 is true if $p(z_1,t_1)$ is replaced with $p(z_1,t_1)$). Let us assume that
$k>2$. We may express $I^k(f_1,\dots,f_k)$ as
\begin{equation*}
\int_0^t\int_0^{t_k}\int_{\bR^{2d}}A(z_{k-1},t_{k-1})f_k(z_k,t_k)
p^{(2)}(z_k-z_{k-1},t_k-t_{k-1})(t-t_k)^\a\ dz_{k-1}dz_k dt_{k-1}dt_k,
\end{equation*}
where $A(z_{k-1},t_{k-1})$ is given by
\begin{equation*}
\int_{\D^{k-2}(t_{k-1})}\int_{\bR^{(k-2)d}} f_1(z_1,t_1)p(z_1,t_1)\prod_{i=2}^{k-1}f_i(z_i,t_i)
p^{(1)}_i(z_i-z_{i-1},t_i-t_{i-1})\ \prod_{i=1}^{k-2}dz_idt_i .
\end{equation*}
By Lemma~4.2 the expression $A(z_{k-1},t_{k-1})$ is bounded above by
\[
\le C_9^{k-2}\ \nu^{-\left(1+\frac dr\right)\frac {k-2} 2}\ f_{k-1}(z_{k-1},t_{k-1})
\ p\left(z_{k-1},2t_{k-1}\right)\left((k-2) \d_1\right)^{-(k-2) \d_1}
t_{k-1}^{(k-2) \d_1}\prod_{i=1}^{k-2}\|f_i\|_{r,q} .
\]
Hence, 
\begin{align*}
I^k(f_1,\dots,f_k)=
\int_0^t\int_0^{t_k}\int_{\bR^{2d}}g(z_{k-1},t_{k-1})&f_k(z_k,t_k)p(z_{k-1},2t_{k-1})
p^{(2)}(z_k-z_{k-1},t_k-t_{k-1})\\
&\x t_{k-1}^{(k-2) \d_1}(t-t_k)^\a\ dz_{k-1}dz_k dt_{k-1}dt_k,
\end{align*}
where $g=f_{k-1}B$ with
\begin{equation}\label{eq4.10}
\left|B(z_{k-1},t_{k-1})\right|\le  C_9^{k-2}\ \nu^{-\left(1+\frac dr\right)\frac {k-2} 2}\
\left((k-2) \d_1\right)^{-(k-2) \d_1}\prod_{i=1}^{k-2}\|f_i\|_{r,q}.
\end{equation}
Since $I^k$ can be written as $I'(g,f_k)$, we can use Lemma~4.1 and \eqref{eq4.10} to assert
that the expression $|I^k(f_1,\dots,f_k)|$ is bounded above by
\[
 C_8 C_9^{k-2}\ \nu^{-\left(1+\frac dr\right)\frac {k} 2}\
\left((k-2) \d_1\right)^{-(k-2) \d_1}\frac{\a^\a \left((k-2) \d_1\right)^{(k-2) \d_1}}
{(\a+k\d_1)^{\a+k\d_1}}\left((k-2) \d_1+1\right)^{ \d_1-\frac d{2r}}\
t^{\a+k \d_1}\prod_{i=1}^{k}\|f_i\|_{r,q}.
\]
From this we can readily deduce \eqref{eq3.2}.
\qed

\subsection{Bounding Double Integrals}

The main reason that we were able to bound the block integrals $I(\b_1,\dots,\b_n)$ that appeared
in \eqref{eq3.3} has to do with the fact that $\b_1+\dots+\b_n=n$. This means that any second derivative of $p$ much be matched with a $0$-th derivative so that the singular integral associated with a second derivative can be controlled. However, if in place of \eqref{eq1.5} we assume the stronger condition \eqref{eq1.17},
then bounding the block integrals of type $I_k$ becomes easier because we can bound double integrals involving first and second derivatives of $p$. To explain this, let us define 
 $K(f_1,f_2)$ as 
\[
\int_{\D^2}\int_{\bR^{2d}} f_1(z_1,t_1)f_2(z_2,t_2)p^{(1)}(z_1,t_1-t_0)
 p^{(2)}(z_2-z_1,t_2-t_1) (t-t_2)^\a\ dz_1dz_2dt_1dt_2,
\]
where $p^{(1)}=p_{z^i}$ and $p^{(2)}=p_{z^jz^k}$ for some $i,j,k\in\{1,\dots,d\}$.

\begin{theorem}\label{thm4.1}Assume \eqref{eq1.17}. There exists a constant $C_{10}=C_{10}(r,q)$ 
such that
\begin{equation}\label{eq4.11}
|K(f_1,f_{2})|\le C_{10}\nu^{-\left(\frac d{2r}+\frac 12\right)}\hat\z(\a)(t-t_0)^{\a+2\d_2}
\|f_1\|_{r,q} \|f_2\|_{r,q},
\end{equation}
where
\[
\hat\z(\a)=\frac{\a^\a}{(\a+{2\d_2})^{\a+2\d_2}}.
\]
\end{theorem}

{\bf{Proof.}} 
The proof is only sketched because it is very similar to the proof of Lemma~4.1 when $k=2$.
Without loss of generality, assume that $t_0=0$, and define $f_i'$ as in the proof of Lemma~4.1.
We decompose $K=K_1+K_2$ where $K_i$ is obtained from $K$ by replacing the domain of integration
 $\D^2$ with 
$\D_i,$ and $\D_1$ and $\D_2$ are defined as in the proof of  Lemma~4.1.
The expression $|K_1|$ is bounded above by
\begin{align*}
 c_0&\nu^{-\left(\frac dr+\frac 32\right)}\int_{\D_1}f_1'(t_1)f_2'(t_2) t_1^{-\frac d{2r}-\frac 12}
(t_2-t_1)^{-\frac d{2r}-1}(t-t_2)^{\a }\ dt_1dt_2\\
&\le c_0\nu^{-\left(\frac dr+\frac 32\right)}\|f_1\|_{r,q}\|f_2\|_{r,q}\ 
\left(\int_{\D_1} t_1^{-\left(\frac {d}{2r}+\frac 12\right)q'}(t_2-t_1)^{-\left(\frac d{2r}+1\right)q'}(t-t_2)^{\a q'}
\ dt_1dt_2\right)^{\frac 1{q'}}\\
&\le c_0\nu^{-\left(\frac dr+\frac 32\right)}\|f_1\|_{r,q}\|f_2\|_{r,q}\ 
\left(\int_{\D_1} t_1^{-\left(\frac d{2r}+\frac 34\right)q'}
(t_2-t_1)^{-\left(\frac d{2r}+\frac 34\right)q'}(t-t_2)^{\a q'}\ dt_1dt_2\right)^{\frac 1{q'}}\\
&=c_0\nu^{-\left(\frac dr+\frac 32\right)}\|f_1\|_{r,q}\|f_2\|_{r,q}\ t^{2\d_1+\a}\eta(\a),
\end{align*}
 where  $\eta$ is defined by
\[
\eta(\a):=\left(\frac {\G(\d_2q')^2\G(\a q' +1)}{\G(2\d_2q'+\a q'+1)}\right)^{\frac 1{q'}},
\]
with $q'=q/(q-1)$. 
As a result,
\begin{equation}\label{eq4.12}
|K_1|\le c_0\nu^{-\left(\frac dr+\frac 32\right)}\eta(\a)\|f_1\|_{r,q}\|f_2\|_{r,q}\ t^{2\d_2+\a},
\end{equation}
It remains to bound $K_2$.

{\bf{Step 2.}} We next decompose $K_2$ as $K_{21}+K_{22}$, where $K_{21}$ and $K_{22}$ 
are defined as in Step 2 of the proof of Lemma~4.1. 
This time the corresponding $I_{2}(k,\ell)$ is defined by
\[
\int_{\D_2}\int_{B_{k\ell}(t_1,t_2)} f_1(z_1,t_1)f_2(z_2,t_2)p^{(1)}(z_1,t_1)
 p^{(2)}(z_2-z_1,t_2-t_1)(t-t_2)^\a\ dz_1dz_2dt_1dt_2.
\]
Again if $(z_1,z_2)\in B_{k\ell}(t_1,t_2)$, with $(t_1,t_2)\in \D_2$ and some $(k,\ell)\in \L_1$, then 
\[
|z_2-z_1|^2\ge  d{\nu t_1}.
\]
As a result,
\begin{align*}
| p^{(2)}(z_2-z_1,t_2-t_1)|&\le \frac {c_1}{\nu(t_2-t_1)}p(z_2-z_1,2(t_2-t_1))\\
&= \frac {c_1}{(\nu(t_2-t_1))^{\frac 34}} 
\frac {|z_2-z_1|^{\frac 12}}{(\nu(t_2-t_1))^{\frac 14}}\frac 1{|z_2-z_1|^{\frac 12}}p(z_2-z_1,2(t_2-t_1))\\
&\le \frac {c_2}{(\nu(t_2-t_1))^{\frac 34}} \frac 1{|z_2-z_1|^{\frac 12}}p(z_2-z_1,4(t_2-t_1))\\
&\le \frac {c_2}{(\nu(t_2-t_1))^{\frac 34}} \frac 1{(\nu t_1)^{\frac 14}}p(z_2-z_1,4(t_2-t_1)).
\end{align*}
This in turn implies that the term $K_{21}$ is bounded above by a constant multiple of 
\begin{align*}
\nu^{-1}&\int_{\D}\int_{\bR^{2d}}| f_1(z_1,t_1)f_2(z_2,t_2)|t_1^{-\frac 14}p^{(1)}(z_1,2t_1)|
 (t_2-t_1)^{-\frac 34}p(z_2-z_1,4(t_2-t_1))\\
&\ \ \ \ \ \ \ \ \ \ \ \ \ \ \ \ \ \ \ \ \ \ \ \ \times (t-t_2)^\a\ dz_1dz_2dt_1dt_2\\
&\le c_3\nu^{-\left(\frac dr+\frac 32\right)}\int_{\D}f_1'(t_1)f_2'(t_2) t_1^{-\frac d{2r}-\frac 34}
(t_2-t_1)^{-\frac d{2r}-\frac 34}(t-t_2)^{\a }\ dt_1dt_2\\
&\le c_3\nu^{-\left(\frac dr+\frac 34\right)}\|f_1\|_{r,q}\|f_2\|_{r,q}\ t^{2\d_2+\a}\eta(\a).
\end{align*}
As a result,
\begin{equation}\label{eq4.13}
|I_{21}|\le c_{4}\nu^{-\left(\frac dr+\frac 34\right)}\eta(\a)
\|f_1\|_{r,q}\|f_2\|_{r,q}\ t^{2\d_2+\a},
\end{equation}
It remains to bound $K_{22}$.

{\bf{Step 3.}} Define $f_{1k}$ and $f_{2\ell}$ as in Step 3 of the proof of Lemma~4.1 and use
  Plancheral's formula, to assert that for any $\d>0$, the term $|K_2(k,\ell)|$ is bounded above by
\begin{align*}
 2 {\pi^2}\int_{\D_2}\int_{\bR^{d}}& \left( \d^{-1}(\nu t_1)^{\frac {d+1}2}|\hat f_{1k}(\xi,t_1)|^2+
\d (\nu t_1)^{-\frac {d+1}2}|\check f_{2\ell}(\xi,t_2)|^2\right)|\xi|^2\\
&\ \ \ \ \ \ \ \ \ \times e^{-4\pi^2\nu(t_2-t_1)|\xi|^2}(t-t_2)^\a\ d\xi  dt_1dt_2
=:K_2^1(k,\ell)+K_2^2(k,\ell).
\end{align*}
Further, $K_2^1(k,\ell)$ bounded above by
\begin{align*}
& \frac {2\pi^2}{\d}\int_{\D}\int_{\bR^{d}}(\nu t_1)^{\frac {d+1}2}|\hat f_{1k}(\xi,t_1)|^2
|\xi|^2e^{-4\pi^2\nu(t_2-t_1)|\xi|^2}(t-t_1)^\a\ d\xi dt_1dt_2\\
&\le \frac {1}{2\nu\d}\int_{0}^t\int_{\bR^{d}}(\nu t_1)^{\frac {d+1}2}|\hat f_{1k}(\xi,t_1)|^2(t-t_1)^\a
\ d\xi dt_1\\
&=\frac {1}{2\nu\d}\int_{0}^t\int_{\bR^{d}}(\nu t_1)^{\frac {d+1}2}| f_{1k}(x,t_1)|^2(t-t_1)^\a
\ dx dt_1\\
&\le  {c_{5}e^{-\frac 14(|k|-\sqrt d)^{+2}}}(\nu \d )^{-1}\int_{0}^t 
(\nu t_1)^{-\left(\frac d2+\frac 12\right)}\int_{B_k(t_1)}| f(x,t_1)|^2(t-t_1)^\a
\ dx dt_1.
\end{align*}
We now apply H\"older's inequality to assert that $K_2^1(k,\ell)$ bounded above by
\begin{equation}\label{eq4.14}
|K_{2}^1(k,\ell)|\le c_{6}\nu^{-\left(\frac dr+\frac 34\right)}
e^{-\frac 15|k|^2}\eta'(\a)\|f_1\|_{r,q}^2\ t^{2\d_2+\a}.
\end{equation}
where 
\[
\eta'(\a)=\left(
\frac{\G(2\d_2q/(q-2))\G(\a q/(q-2)+1)}{\G((2\d_2+\a)q/(q-2)+1)}\right)^{\frac {q-2}q}.
\]
In the same fashion we show
\begin{equation}\label{eq4.15}
|K_{2}^2(k,\ell)|\le c_{7}{\d }\nu^{-\left(\frac dr+1\right)}\eta'(\a)\|f_2\|_{r,q}^2\ t^{2\d_2+\a}.
\end{equation}
The rest of the proof is as in the proof of Lemma~4.1.
\qed

\section{Symplectic Diffusions and Navier-Stokes Equation}
\label{sec5}
 \setcounter{equation}{0}

\noindent
{\bf{Proof of Theorem~1.2. Step 1.}}  For $u\in L^{r,q}$, with $r$ and $q$ satisfying \eqref{eq1.5}, choose a sequence of smooth functions $u_N$ such that 
$\|u_N-u\|_{r,q}\to 0$ as $N\to\i$.
Write $\O$ for the space of pair of continuous functions $X:\bR^d\x[0,T]\to\bR^d$
and  $B:[0,T]\to\bR^d$ such that  $X$ is weakly differentiable with respect 
to the spatial variables and  $D_aX$ is
 locally in $L^p$ for every $p\ge 1$. We equip $\O$ with a topology of $L^\i_{loc}$ for $X$ and weak topology for $D_aX$. Consider the SDE 
\begin{equation}\label{eq5.1}
dX=u_N(X,t)dt+\s dB,
\end{equation}
where $B$ is a standard Brownian motion. The law of the pair $(X,B)$ 
is a probability measure $\cP^N$ on the space $\O$ such that the $B$ 
component is a standard Brownian motion. 
 Using the equations \eqref{eq5.1}, \eqref{eq2.2}, \eqref{eq2.10} and Girsanov's formula we can readily show
\begin{equation}\label{eq5.2}
\int\sup_{t\in[0,T]}|X(a,t)-a|^2\ d\cP^N\le c_0T+c_0\int\int_0^T|u_N(X(a,t)|^2 dt\ d\cP^N\le c_1T,
\end{equation}
for a constant $c_1$ independent of $N$. We may use \eqref{eq5.2}, 
Theorem~1.1 and Corollary~1.2 to assert that the family $\{\cP_N\}_{N=1}^\i$ is tight.
Let $\cP$ be a limit point of the family $\{\cP_N\}_{N=1}^\i$ as $N\to\i$. 
Let $J:\bR^d\to [0,\i)$ be a 
continuous function of compact support. Use \eqref{eq5.1} to assert
\begin{align*}
\lim_{N\to\i}&\int\left\{\sup_{t\in[0,t]}\int_{\bR^d}
\left|X(a,t)-a-B(t)-\int_0^t u(X(a,s),s)\ ds\right|\ J(a)\ da\right\}\ d\cP_N\\
&=\lim_{N\to\i}
\int\left[\sup_{t\in[0,t]}\int_{\bR^d}\left|\int_0^t (u_N-u)(X(a,s),s)\ ds\right|\ J(a)\ da\right]\ d\cP_N\\
&=\lim_{N\to\i}
\int\left[\sup_{t\in[0,t]}\int_{\bR^d}\left|\int_0^t (u_N-u)(a,s)\ ds\right|\ J\left(X^{-1}(a,t)\right)
\left|\det D_aX^{-1}(a,t)\right|\ da\right]\ d\cP_N=0.
\end{align*}
Note that the expression inside the curly brackets is a continuous functional. As a result, we may use
our bounds on $D_aX$ to show 
\[
\int\left\{\int_{\bR^d}\left|X(a,t)-a-B(t)-\int_0^t u(X(a,s),s)\ ds\right|\ J(a)\ da\right\}\ d\cP
=0,
\]
which in turn implies that the equation 
\[
X(a,t)=a+\int_0^t u(X(a,s),s)\ ds+B(t),
\]
 is valid $\cP$-almost surely for almost all $a\in\bR^d$, and hence for all $a$ by continuity.

{\bf{Step 2.}} We now verify \eqref{eq1.10}. Since $u_N$ is smooth, we apply Proposition~3.1 of [R]
to assert 
\begin{align}\nonumber
\int\left(X^*_t\b^t-\b^0\right)(V)\ dx&-\int_0^t 
\left(
\int X^*_s\left[\dot\b^s+\cA_{u_N}\b^s\right](V)\ dx \right)\ ds\\
&\ \ \ \ \ \ \ =\int_0^t \int\left(\sum_iX^*_s\g_i^s\right)(V)\ dx\ dB^i(s),
\label{eq5.3}
\end{align}
$\cP_N$-almost surely. Replacing $u_N$ with $u$ results in an error that is bounded above by a 
constant multiple of 
\begin{align*}
Err_t(X):=&\int_0^t\int|D_aX(a,s)|\ |(u_N-u)(X(a,s),s)|\ |V(a)|\ dads\\
&+\int_0^t\int |(u_N-u)(X(a,s),s)|\ |\nabla\cdot V(a)|\ dads.
\end{align*}
Finally we use \eqref{eq1.6} and \eqref{eq5.2} to show
\[
\lim_{N\to\i}\int \sup_{t\in[0,T]}Err_t\ d\cP_N=0.
\]
This allows us to pass to the limit in \eqref{eq5.3} and deduce \eqref{eq1.10}.
\qed

\bigskip
\noindent
{\bf{Proof of Corollary~1.3}} Let us write $x=(q,p)$ and set $\l=p\cdot dq$. We certainly have
\[
\cA_u\l=-\hat d(H-p\cdot H_p),\ \ \ w_{q^i}=0,\ \ \ w_{p^i}\cdot dx=dp^i,
\]
where $w=[p,0]$. As a result the forms $\cA_u\l$ and $w_{x^i}\cdot dx$ are exact for $i=1,\dots,n$. From this and
Theorem~1.2 we learn that $X_t^*\l$ is exact. This in turn implies that $X_t^*\hat d\l=0$, as desired.
\qed

Given a classical solution $u(\cdot,t)$ of \eqref{eq1.1}, let us write $\a^t=u(\cdot,t)\cdot dx$
for the $1$-form associated with $u$. In terms of $\a$, the equation \eqref{eq1.1} may be written as
\begin{equation}\label{eq5.4}
\dot\a^t+i_u\hat d\a^t+\hat dH^t=0,
\end{equation}
where $H^t(x)=\frac 12|u(x,t)|^2+P(x,t)$. Here $i_u$ denotes the contraction operator and
we are simply using the identity 
\[
\sum_ju^i_{x^j}u^j=\sum_j(u^i_{x^j}-u^j_{x^i})u^j+\left(\frac 12|u|^2\right)_{x^i}.
\]
Further, if we use Cartan's formula and write $\cL_u=\hat d\circ i_u+i_u\circ \hat d$ for the Lie derivative in the direction of u, we may rewrite \eqref{eq5.4} as
\begin{equation}\label{eq5.5}
\dot\a+\cL_u\a^t-\hat dL^t=0,
\end{equation}
where $L=\frac 12|u|^2-P.$  Equation \eqref{eq5.5} can be used to give a geometric description of the Euler Equation \eqref{eq1.1}: If we write  $X(\cdot,t)=X_t(\cdot)$ for the flow 
of $u$ as in \eqref{eq1.2}, then \eqref{eq5.3} really means 
\[
\frac d{dt}X^*_t\a^t=X^*_t\hat d L^t,
\]
or equivalently
\begin{equation}\label{eq5.6}
X_t^*\a^t-\a^0=\hat dK^t,
\end{equation}
for $K^t=\int_0^t L^s\circ X_s\ ds$. The identity \eqref{eq5.6} is the celebrated Kelvin's circulation formula
and coupled with the incompressibility condition $\nabla\cdot u=0$ is equivalent to Euler Equation. 

In the case of viscid fluid, the fluid velocity satisfies Navier-Stokes Equation. For our purposes, it is more convenient to use backward Navier-Stokes Equation
\begin{equation}\label{eq5.7}
u_t+(Du)u+\nabla P+\nu\D u=0.
\end{equation}
For a classical solution of \eqref{eq5.7}, we may write
\begin{equation}\label{eq5.8}
\dot\a^t+\cA_u\a^t-\hat dL=0.
\end{equation}
On the other hand, if $X_t$ denotes the flow of SDE \eqref{eq1.4} and 
$\b^t=X^*_t\a^t$, then
\[
M^t:=\b^t-\b^0-\int_0^t X^*_s(\dot\a+\cA_u\a)\ ds=\b^t-\b^0-\hat dK^t,
\] 
is a martingale.
In summary 
\begin{equation}\label{eq5.9}
X^*_t\a^t=\a^0+M^t+\hat dK^t.
\end{equation}
By taking the exterior derivative, we obtain
\begin{equation}\label{eq5.10}
X^*_t\hat d\a^t=\hat d\a^0+\hat d M^t.
\end{equation}
For both \eqref{eq5.9} and \eqref{eq5.10} we are assuming that $u$ is a classical solution of Navier-Stokes Equation. For a weak solution of \eqref{eq5.7} we wish to show that $M^t$ is still a martingale.
This is exactly the content of Theorem~1.3 provided that the solution $u$ can be approximated by suitable regular functions.

\bigskip
\noindent
{\bf{Proof of Theorem~1.3.}} Assume that $(v,w)=(v^\e,w^\e)$ solves Camassa-Holm Equation
with $v=w-\e\D w$. Set $\bar\a^t=v(\cdot,t)\cdot dx$ and write $Y$ for the flow of the SDE
\[
dY=w(Y,t)dt+\s dB.
\]
As in \eqref{eq5.6},
 the equation \eqref{eq1.14} can be rewritten as
\begin{equation*}
\frac{d}{dt}\bar\a^t+\cA_w\bar\a^t-\hat d\bar L^t=0.
\end{equation*}
This in turn implies 
\begin{equation}\label{eq5.11}
Y^*_t\bar\a^t=\bar\a^0+\bar M^t+\hat d\bar K^t,\ \ \ Y^*_t\hat d\bar\a^t=\hat d\bar\a^0+\hat d
\bar M^t,
\end{equation}
where $\bar M^t$ is a martingale and
$\bar K^t=\int_0^t \bar L^s\circ Y_s\ ds$. We now choose a subsequence of 
$w=w^\e$ so that $w^\e\to u$. From our assumption
 \eqref{eq1.15}, Theorem~1.1 and Corollary~1.2 we can choose a further subsequence such that  $Y=Y^\e\to X$ in $L^\i_{loc}$, and $D_aY^\e\to D_aX$
weakly in any $L^p$ space. This allows us to pass to the limit in \eqref{eq5.11} to assert that
the process \eqref{eq1.16} is a martingale.
\qed

\end{document}